\documentclass[twoside,a4paper,11pt,leqno]{amsart}

\usepackage{amssymb}
\swapnumbers
\theoremstyle{plain}
\newtheorem{thm}{Theorem}[section]
\newtheorem{prop}[thm]{Proposition}

\newtheorem{lemma}[thm]{Lemma}
\theoremstyle{definition}
\newtheorem{defn}[thm]{Definition}
\theoremstyle{remark}
\newtheorem{rem}[thm]{Remark}
\newtheorem{examples}[thm]{Examples}
\newcommand{\proofof}[1]{\end{#1}\begin{proof}}
\newcommand{\emphdef}{\textit}
\numberwithin{equation}{section}

\newcommand{\acknowledge}{\section*{Acknowledgements}}
\newcommand{\bauth}[1]{\textit{#1}, }
\newcommand{\bart}[1]{#1, }
\newcommand{\bjourn}[3]{#1 \textbf{#2} (#3) }
\newcommand{\bbook}[1]{#1, }

\newcommand{\bpp}[2]{#1--#2.}
\newcommand{\tcite}[1]{\textup{\cite{#1}}}
\DeclareMathAlphabet{\mathrmsl}{OT1}{cmr}{m}{sl}
\newcommand{\rssymb}[2]{\newcommand{#1}{{\mathrmsl{#2\mkern1mu}}}}
\newcommand{\calsymb}[2]{\newcommand{#1}{{\mathcal{#2}}}}
\newcommand{\bbsymb}[2]{\newcommand{#1}{{\mathbb{#2}}}}
\newcommand{\lieoper}[2]{\newcommand{#1}{\mathop{\mathfrak{#2}}}}
\newcommand{\oper}[3][n]{\newcommand{#2}{\mathop{\mathrm{#3}}\ifx
  n#1\nolimits\else\limits\fi}}
\newcommand{\rsoper}[3][n]{\newcommand{#2}{\mathop{\mathrmsl{#3\mkern1mu}}\ifx
  n#1\nolimits\else\limits\fi}}
\bbsymb\C{C}\bbsymb\HQ{H}\bbsymb\R{R}
\calsymb\cD{D}\calsymb\cG{G}\calsymb\cH{H}\calsymb\cK{K}\calsymb\cL{L}
\calsymb\cN{N}\calsymb\cW{W}
\newcommand{\eps}{\varepsilon}
\newcommand{\gam}{\gamma}
\newcommand{\lam}{\lambda}
\oper\End{End}\oper\Sym{Sym}\oper\Cl{Cl}
\oper\Aut{Aut}\oper\GL{GL}\oper\SL{SL}
\oper\CO{CO}\oper\On{O}\oper\SO{SO}
\oper\CSpin{CSpin}\oper\Symp{Sp}\oper\Un{U}\oper\SU{SU}\oper\CU{CU}
\lieoper\gl{gl}\lieoper\sgl{sl}
\lieoper\co{co}\lieoper\on{o}\lieoper\so{so}
\lieoper\symp{sp}\lieoper\un{u}\lieoper\su{su}\lieoper\csu{csu}
\newcommand{\ip}[1]{\langle#1\rangle}
\newcommand{\Ip}[1]{\bigl\langle#1\bigr\rangle}
\newcommand{\lie}[1]{\mathfrak{#1}}
\renewcommand{\geq}{\geqslant}
\renewcommand{\leq}{\leqslant}
\newcommand{\dsum}{\oplus}               
\newcommand{\Dsum}{\bigoplus}            
\newcommand{\tens}{\mathbin{\otimes}}    
\newcommand{\Tens}{\bigotimes}           
\newcommand{\cartan}{\mathbin{\odot}}    
\newcommand{\idealsub}{\rtimes}          
\newcommand{\subnormal}{\ltimes}         
\newcommand{\normalsub}{\rtimes}         
\newcommand{\intersect}{\mathinner{\cap}}
\newcommand{\act}{\mathinner{\cdot}}
\newcommand{\dual}{^{*\!}}
\newcommand{\Cinf}{\mathrm{C}^\infty}
\rsoper\Ad{Ad}
\rsoper\kernel{ker}
\rsoper\image{im}
\rsoper\alt{alt}           
\rsoper\sym{sym}           
\rsoper\trace{tr}          
\rsoper\detm{det}          
\rsoper\divg{div}
\rsoper\Sdivg{Sdiv}
\rsoper\Twist{Twist}
\rsoper\project{proj}
\rsoper\represent{repr}
\rssymb\ev{ev}
\rssymb\ad{ad}
\rssymb\iden{id}
\newcommand{\cross}{\mathbin{{\times}\!}\low}
\newcommand{\from}{\colon}
\newcommand{\isom}{\cong}
\newcommand{\hQuabla}{\hat{\pmb\square}}
\newcommand{\Quabla}{\pmb{\square}}

\newcommand{\conf}{\mathsf{c}}

\newcommand{\cform}{\eta}
\newcommand{\g}{\lie{g}}
\newcommand{\p}{\lie{p}}
\newcommand{\Lieb}[1]{[#1]}
\bbsymb\E{E}\bbsymb\W{W}
\newcommand{\T}{\lie{m}}
\newcommand{\inner}{\mathbin{\lrcorner}}
\newcommand{\capinner}{\mathbin{\lower3pt\hbox{$\urcorner$}}}
\newcommand{\low}{^{\vphantom x}}
\newcommand{\dT}{d\low_{\T}}
\newcommand{\dTd}{d\low_{\T^*}}
\newcommand{\delT}{\delta\low_{\T}}
\newcommand{\delTd}{\delta\low_{\T^*}}
\newcommand{\delTdM}{\delta\low_{T\dual M}}
\newcommand{\gM}{\g\low_M}
\newcommand{\InvDer}{\nabla^\cform}
\newcommand{\TwConn}{\nabla^\g}
\newcommand{\deltw}{\delta^\g}
\newcommand{\dtw}{d^\g}
\newcommand{\Rtw}{R^\g}
\newcommand{\QuaTw}{\Quabla_\g}
\newcommand{\Qtw}{Q}
\newcommand{\PiTw}{\Pi}
\newcommand{\OpTw}{\cD}
\newcommand{\cptw}{\mathbin{\scriptstyle\sqcup}}
\newcommand{\delinv}{\delta^\cform}
\newcommand{\dinv}{d^\cform}
\newcommand{\Rinv}{R^\cform}
\newcommand{\hQuaInv}{\hQuabla_\cform}
\newcommand{\QuaInv}{\Quabla_\cform}
\newcommand{\Qinv}{Q_\cform}
\newcommand{\PiInv}{\Pi^\cform}
\newcommand{\OpInv}{\cD^\cform}
\newcommand{\Opinv}{\cD_\cform}
\newcommand{\cpinv}{\cptw\low_\cform}
\newcommand{\TwIso}{\Psi_\W}
\newcommand{\captw}{\mathbin{\scriptstyle\sqcap}}
\newcommand{\capinv}{\captw\low_\cform}
\newcommand{\Dirac}{{\smash{\raise1pt\hbox{$\not$}}\mkern-1.5mu D}}
\begin{document}
\title[Differential invariants and curved BGG sequences]{Differential
invariants and curved Bernstein-Gelfand-Gelfand sequences}
\author{David M. J. Calderbank}
\address{Department of Mathematics and Statistics\\ University of Edinburgh\\
King's Buildings, Mayfield Road\\ Edinburgh EH9 3JZ. Scotland.}
\email{davidmjc@maths.ed.ac.uk}
\author{Tammo Diemer}
\address{Mathematisches Institut, Universit\"at Bonn\\
Beringstra\ss e 1\\ D-53113 Bonn. Germany.}
\curraddr{DePfa Deutsche Pfandbrief Bank AG\\ Paulinenstra\ss e 15\\
D-65189 Wiesbaden. Germany.}
\email{tammo.diemer@depfa.com}
\date{December 1999; revised February 2000}
\keywords{Bernstein-Gelfand-Gelfand resolution, Cartan connection,
parabolic geometry, cup product, Leibniz rule, strong homotopy algebra,
twistor methods}
\begin{abstract}
We give a simple construction of the Bernstein-Gelfand-Gelfand sequences of
natural differential operators on a manifold equipped with a parabolic
geometry. This method permits us to define the additional structure of a
bilinear differential ``cup product'' on this sequence, satisfying a Leibniz
rule up to curvature terms. It is not associative, but is part of an
$A_\infty$-algebra of multilinear differential operators, which we also obtain
explicitly.  We illustrate the construction in the case of conformal
differential geometry, where the cup product provides a wide-reaching
generalization of helicity raising and lowering for conformally invariant
field equations.
\end{abstract}
\maketitle
\vspace{-4mm}
\section*{Introduction}

In a sequence of pioneering papers~\cite{Baston1,Baston2,Baston3}, Robert
Baston introduced a number of general methods to study invariant differential
operators on conformal manifolds, and a related class of parabolic geometries,
which he called ``almost hermitian symmetric (AHS) structures''. In
particular, he suggested that certain complexes of natural differential
operators, dual to generalized Bernstein-Gelfand-Gelfand (BGG) resolutions of
parabolic Verma modules, could be extended from the homogeneous context
(generalized flag manifolds) to curved manifolds modelled on these spaces. He
provided a construction of such a BGG sequence (no longer a complex in
general) for AHS structures~\cite{Baston2}, and introduced (in~\cite{Baston3})
a class of induced modules, now called semiholonomic Verma modules~\cite{ES}.

Baston's work fits into the programme of parabolic invariant theory initiated
by Fefferman and Graham~\cite{Fef,FG}. Several authors have joined in an
endeavour to complete these ideas and hence provide a theory of invariant
operators in all parabolic geometries, which include conformal geometry,
projective geometry, quaternionic geometry, projective contact geometry, CR
geometry and quaternionic CR geometry.  In~\cite{ES}, Eastwood and Slov\'ak
began the study of semiholonomic Verma modules and classified the Verma module
homomorphisms lifting to the semiholonomic modules in the conformal case. The
AHS structures have been extensively studied by \v Cap, Slov\'ak and Sou\v cek
in~\cite{CSS1,CSS2,CSS3}, and in the last paper of this series they construct
a large class of invariant differential operators for these geometries.
Then, in~\cite{CSS4}, \v Cap, Slov\'ak and Sou\v cek clarified
Baston's construction of the BGG sequences in the AHS case, and in the
process, generalized it to all parabolic geometries.  Hence we now know that
all standard homomorphisms of parabolic Verma modules induce differential
operators also in the curved setting, providing us with a huge supply of
invariant linear differential operators.

This paper has two main objectives: to simplify the construction of the BGG
sequences given in~\cite{CSS4}, and to equip these sequences of linear
differential operators with bilinear differential pairings, inducing, in the
flat case, a cup product on cohomology.

The key tool for the construction of the BGG sequences is an invariant
differential operator from relative Lie algebra homology bundles to twisted
differential forms, denoted $L$ in~\cite{Baston2} and~\cite{CSS4}.  However,
when one tries to produce a cup product, one needs a differential operator
defined on the whole bundle of twisted differential forms which induces $L$ on
the homology bundles. The search for such an operator (within the dual
homogeneous formalism of Verma modules) led the second author to a procedure
which, in addition to providing a cup coproduct on the BGG resolutions, gives
a simpler construction of the resolutions themselves.  These developments are
described in~\cite{Thesis}.  It is straightforward to dualize this
procedure and one readily sees that it generalizes from the homogeneous
context to arbitrary parabolic geometries, although the presence of curvature
introduces phenomena that do not arise in the flat case, and also suggests
further constructions of multilinear differential operators.  We present, in
geometric language, these constructions and phenomena here.  That is, we give
a simple self-contained approach to the curved BGG sequences and cup
product, in their natural geometrical context, and introduce an
$A_\infty$-algebra of invariant multilinear differential operators.

We begin by recalling some basic facts from Cartan geometry, emphasizing the
simple first order constructions a Cartan connection provides.  Our approach
is mainly influenced by~\cite{Baston1,CSS4,Sharpe,Slovak}.  In
section~\ref{pg}, we define parabolic geometries as Cartan geometries
associated to a semisimple Lie algebra $\g$ with a parabolic Lie subalgebra
$\p$. We summarize the basic representation-theoretic facts that we will need
and give some examples.

The most substantial piece of representation theory we need is Lie algebra
homology, and we discuss this in section~\ref{lahc}. In order to keep the
paper as self-contained as possible, we give proofs for all the basic
properties of Lie algebra homology, although we only indicate briefly how
Kostant's Hodge decomposition may be established. Additionally, there are some
non-standard aspects to our treatment: firstly we concentrate on Lie algebra
homology, rather than cohomology, since it is the Lie algebra homology that is
$\p$-equivariant, and secondly, we eschew the lamentably inverted notation
$\partial,\partial^*$ for the Lie algebra coboundary and boundary operators
(for some reason, $\partial$, although a boundary operator in~\cite{Kostant},
is the coboundary operator in~\cite{Baston1,CS,CSS2,Tanaka}).  Instead,
following Kostant in part~\cite{Kostant}, and by analogy with the deRham
complex, we use $d$ and $\delta$, with subscripts to indicate that it is the
Lie algebraic rather than differential operators we are considering. This
analogy with the deRham complex is central to our proof. After stating the
main theorem to be proved at the end of section~\ref{lahc}, we begin the
study, in section~\ref{tdr}, of the twisted deRham complex. As observed
in~\cite{CSS4}, there are in fact two natural deRham complexes one might
consider, which differ if the parabolic geometry has torsion.

We prove an explicit version of our main result in section~\ref{bgg}.  There
we give a construction, using a Neumann series, of the BGG sequences of
differential operators found in~\cite{CSS4}, and at the same time construct
the bilinear differential pairings.  Our method enables us to compute
explicitly the curvature terms entering into the squares of the BGG
differentials and into the Leibniz rule for the pairings. The BGG sequence
of~\cite{CSS4} is based on the the twisted deRham complex with torsion. We
show that under a weak additional assumption, there is another BGG sequence
based on the torsion-free twisted deRham sequence. The operators involved are
in some ways more complicated because of the corrections needed to ``remove''
the torsion, but we believe they are more natural and illustrate this by
interpreting curvature terms for normal regular parabolic geometries. For
torsion-free parabolic geometries, of course, the two BGG sequences agree.

At the end of section~\ref{bgg} and in the following section, we introduce
multilinear differential operators and establish that they form a (curved)
$A_\infty$-algebra.  We study adjointness properties of the BGG operators and
cup product in section~\ref{dualBGG}, introducing a dual BGG sequence and a
cap product.  In section~\ref{genapp}, potential applications, such as
deformation and moduli space problems, are discussed, mostly in a rather
open-ended fashion, since working out the details in many cases is a
substantial research project.  We attempt to be more concrete in the final
section, where we give examples in conformal geometry, and show how the cup
product generalizes helicity raising and lowering by Penrose twistors in four
dimensional conformal geometry to arbitrary twistors in arbitrary dimensions.

Finally, one feature of our methods is that we work with spaces of smooth
sections, and do not need the machinery of semiholonomic infinite jets and
Verma modules. However, for the convenience of the reader, we sketch
an approach to this machine in an appendix.

\acknowledge We are very grateful to Andreas \v Cap, Rod Gover, Martin Markl,
Elmer Rees, Michael Singer, Jan Slov\'ak and Jim Stasheff for stimulating and
helpful discussions. The first author would particularly like to thank
Vladimir Sou\v cek, with whom he has discussed helicity raising and lowering
extensively over the past few years, and more recently, some of the details of
curved BGG sequences.  The second author is similarly indebted to Gregor
Weingart for many discussions on Lie algebra homology and BGG resolutions.

\section{Cartan geometries and invariant differentiation}

Cartan geometries are geometries modelled on a homogeneous space $G/P$ (for a
modern introduction, see~\cite{Sharpe}). Such a homogeneous space has a natural
principal $P$-bundle $G\to G/P$, equipped with $\g$-valued $1$-form
$\omega\from TG\to\g$, namely the Maurer-Cartan form which trivializes $TG$ by
the right-invariant vector fields.

In order to avoid fixing $G$, we work with a pair $(\g,P)$ consisting of a Lie
algebra $\g$ and a group $P$ acting on $\g$ by automorphisms such that the Lie
algebra $\p$ of $P$ is a subalgebra of $\g$ and the action of $P$ on $\g$
induces the adjoint action of $P$ on $\p$ and of $\p$ on $\g$.

\begin{defn}[\textsl{Cartan geometry}] Let $M$ be a manifold of the same
dimension as $\g/\p$.
\begin{enumerate}
\item A \emphdef{Cartan geometry} of type $(\g,P)$ on $M$ is a principal
$P$-bundle $\pi\from \cG\to M$, together with a $P$-equivariant
$\g$-valued $1$-form $\cform\from T\cG\to\g$ such that for each
$u\in\cG$, $\cform_u\from T_u\cG\to\g$ is an isomorphism restricting to
the canonical isomorphism between $T_u(\cG_{\pi(u)})$ and $\p$.
\item A \emph{Kleinian} or \emph{homogeneous model} of a Cartan geometry of
type $(\g,P)$ is a homogeneous space $G/P$, for a Lie group $G$ with
subgroup $P$ and Lie algebra $\g$.
\item The \emphdef{curvature} $K\from\Lambda^2T\cG\to\g$ of a Cartan
geometry is defined by
\begin{equation*}
K(U,V)=d\cform(U,V)+[\cform(U),\cform(V)].
\end{equation*}
It induces a curvature function
$\kappa\from\cG\to\Lambda^2\g\dual\tens\g$ via
\begin{equation*}
\kappa(u)(\xi,\chi)=K_u\bigl(\cform^{-1}(\xi),\cform^{-1}(\chi)\bigr)=
[\xi,\chi]-\cform_u[\cform^{-1}(\xi),\cform^{-1}(\chi)],
\end{equation*}
where $u\in\cG$ and the latter bracket is the Lie bracket of vector fields on
$\cG$.
\item Associated to a $P$-module $\E$ is a vector bundle $E=\cG\cross_P\E$. In
particular, the Cartan connection $\cform$ identifies the tangent bundle of
$M$ with $\cG\cross_P\g/\p$. The \emphdef{adjoint bundle} of a Cartan geometry
is $\gM=\cG\cross_P\g$. The quotient map $\g\to\g/\p$ induces a surjective
bundle map $\pi_{\g}\from\gM\to TM$.

Note that the associated bundle construction identifies sections of
$E=\cG\cross_P\E$ over $M$ with $P$-equivariant maps from $\cG$ to $\E$:
\begin{equation*}
\Cinf(M,E)=\Cinf(\cG,\E)^P.
\end{equation*}
We shall make frequent use of this identification, often without comment.
\end{enumerate}
\end{defn}

The curvature $K$ of a Cartan geometry measures the failure of the
Cartan $1$-form $\cform$ to induce a Lie algebra homomorphism. This is
the obstruction to finding a local isomorphism between $\cG\to M$ and
a homogeneous model $G\to G/P$.

The following definition is essentially given in~\cite{CSS1,CSS4}, except that
we do not restrict the derivative to horizontal tangent vectors, and hence we
do not lose $P$-equivariance. The same idea appears in~\cite{CG,Sharpe}, the
latter reference attributing it to Cartan~\cite{Cartan}.

\begin{defn}[\textsl{Invariant derivative}] Let $(\cG,\cform)$ be a Cartan
geometry of type $(\g,P)$ on $M$, and let $\E$ be a $P$-module with associated
vector bundle $E=\cG\cross_P \E$. Then the \emphdef{invariant derivative}
on $E$ is defined by
\begin{align*}
\InvDer\from
\Cinf(\cG,\E)&\to\Cinf(\cG,\g\dual\tens\E)\\
\InvDer_\xi f &= df\bigl(\cform^{-1}(\xi)\bigr)
\end{align*}
for all $\xi$ in $\g$. It is $P$-equivariant and so maps $\Cinf(\cG,\E)^P$
into $\Cinf(\cG,\g\dual\tens\E)^P$. Identifying $P$-equivariant maps to a
$P$-module with sections of the associated vector bundle therefore provides a
linear map $\InvDer\from \Cinf(M,E)\to\Cinf(M,\g_M\dual\tens E)$.
\end{defn}

We now build up some simple properties of this operation.

\begin{prop}[\textsl{1-jets}]\label{basic} Let $(\cG,\cform)$
be a Cartan geometry of type $(\g,P)$ on $M$.
\begin{enumerate}
\item The curvature $K$ is a horizontal $2$-form and so induces
$K\low_M\in\Cinf(M,\Lambda^2T\dual M\tens\gM)$.  Equivalently $\kappa(\xi,.)=0$
for $\xi\in\p$, so $\kappa\in\Cinf(\cG,\Lambda^2(\g/\p)\dual\tens\g)^P$.
\item The invariant derivative of a $P$-equivariant map $f\from\cG\to\E$ is
vertically trivial in the sense that $(\InvDer_\xi f)(u)
+\xi\act\bigl(f(u)\bigr)=0$ for all $\xi\in\p$ and $u\in\cG$. In particular if
$P$ acts trivially on $\E$ and $f\from M\to E$, then $\InvDer f
=df\circ\pi_{\g}$.
\item If $f_1\from\cG\to \E_1$ and $f_2\from\cG\to \E_2$
then $\InvDer_\xi(f_1\tens f_2)=(\InvDer_\xi f_1)\tens f_2+
f_1\tens(\InvDer_\xi f_2)$.
\item The invariant derivative satisfies the ``Ricci identity'':
\begin{equation*}
\InvDer_\xi(\InvDer_\chi f)-\InvDer_\chi(\InvDer_\xi f)
= \InvDer_{[\xi,\chi]}f - \InvDer_{\kappa(\xi,\chi)}f.
\end{equation*}
\item The map $\,\Cinf(M,E)\to\Cinf\bigl(M,E\dsum(\g_M\dual\tens E)\bigr)$
sending $s$ to $(s,\InvDer s)$ defines an injective bundle map from
$J^1E$ to $E\dsum \bigl(\g_M\dual\tens E\bigr)$ whose image is $\cG\cross_P
J^1_0\E$ where
$J^1_0\E=\{(\phi_0,\phi_1)\in\E\dsum\bigl(\g\dual\tens
\E\bigr):\phi_1(\xi)+\xi\act\phi_0=0 \textup{ for all } \xi\in\p\}$.
\end{enumerate}
\proofof{prop} These are straightforward calculations.
\begin{enumerate}
\item For $\xi\in\p$, we have by definition that $\cform^{-1}(\xi)$ is a
vertical vector field generated by the right $P$-action. Now the map
$\chi\mapsto\cform^{-1}(\chi)$ is $P$-equivariant for any $\chi\in\g$, from
which it follows by differentiating that
$[\cform^{-1}(\xi),\cform^{-1}(\chi)]=\cform^{-1}[\xi,\chi]$.
\item Differentiate the $P$-equivariance condition $p\act(f(up))=f(u)$.
\item This is just the product rule for $d(f_1\tens f_2)$.
\item The Ricci identity holds because both sides are equal to
$df([\cform^{-1}(\xi),\cform^{-1}(\chi)])$.
\item Certainly the map on smooth sections only depends on the $1$-jet at each
point, and it is injective since the symbol of $\InvDer$ is the
inclusion $T\dual M\tens E\to \g_M\dual\tens E$, as one easily sees from the
product rule (for $\E_1$ trivial and $\E_2=\E$). It maps into $\cG\cross_P
J^1_0\E$ by vertical triviality, so the result follows by comparing the ranks
of the bundles.\qed
\end{enumerate}
\begingroup\let\endproof\relax\end{proof}\endgroup
\noindent The final term in the Ricci identity is first order in general, due
to the presence of torsion. The \emphdef{torsion} is defined to be $TM$-valued
$2$-form $\pi_\g(K\low_M)$ obtained by projecting the curvature $K$ onto
$\g/\p$ and a Cartan geometry is said to be \emphdef{torsion-free} if $K$
takes values in $\p$ so that $\pi_\g(K\low_M)=0$ and
$\kappa\in\Cinf(\cG,\Lambda^2(\g/\p)\dual\tens\p)$. In this case, for any
$P$-equivariant $f\from\cG\to\E$, we have
$-\InvDer_{\kappa(\xi,\chi)}f=\kappa(\xi,\chi)\act f\,$.

We end this section by considering the invariant derivative when the
$P$-module is also a $\g$-module.
\begin{defn} A \emphdef{$(\g,P)$-module} is a vector space $\W$ carrying a
representation of $P$ and a $P$-equivariant representation of $\g$, such that
the induced actions of $\p$ coincide.
\end{defn}
\noindent For a Lie group $G$ with Lie algebra $\g$ and subgroup $P$, any
$G$-module is a $(\g,P)$-module.
\begin{defn}[\textsl{Twistor connection}]\label{ltt} Let $\W$ be
a $(\g,P)$-module and define
\begin{align*}
\TwConn\from
\Cinf(\cG,\W)&\to\Cinf(\cG,\g\dual\tens\W)\\
\TwConn_\xi f &= \InvDer_\xi f+\xi\act f\,.
\end{align*}
Then for $P$-equivariant $f$, $\TwConn_\xi f$ vanishes for $\xi\in\p$,
so $\TwConn f$ takes values in $(\g/\p)\dual\tens\W$ and hence
$\TwConn$ induces a covariant derivative on $W=\cG\cross_P\W$ which
will be called the \emphdef{twistor connection} on $W$. Its curvature
is easily computed to be the action of $K\low_M$ on $W$.
\end{defn}
If $G$ is a Lie group with subgroup $P$ and Lie algebra $\g$, then the
principal $G$-bundle $\widetilde\cG=\cG\cross_P G$ has a principal bundle
connection induced by $\cform$, and, on a $G$-module $\W$, $\TwConn$ is simply
the covariant derivative induced by this connection. Although these basic
ideas from the theory of Cartan connections are well-established, the
systematic use of a \emph{linear} representation of the Cartan connection
seems to first appear in twistor theory, where $\W$ is the \emphdef{local
twistor bundle} and $\TwConn$ defines \emphdef{local twistor
transport}. Following Baston~\cite{Baston1}, we adapt this terminology to more
general situations.
\begin{prop}
Let $\TwIso$ the $P$-equivariant automorphism of
$(\E\tens\W)\dsum(\g\dual\tens\E\tens\W)$, for any $P$-module $\E$, sending
$(\phi_0=e\tens w,\phi_1)$ to $(\phi_0,\widetilde\phi_1)$ where
$\widetilde\phi_1(\chi)=\phi_1(\chi)+e\tens\chi\act w$ for any
$\chi\in\g$. Then $\TwIso$ restricts to an isomorphism from $J^1_0(\E\tens\W)$
to $J^1_0(\E)\tens\W$.
\proofof{prop} For any $\chi\in\p$ and $(\phi_0=e\tens w,\phi_1)\in
J^1_0(\E\tens\W)$, we have
\begin{equation*}
\widetilde\phi_1(\chi)+(\chi\act e)\tens w
=\phi_1(\chi)+(\chi\act e)\tens w+e\tens\chi\act w
=\phi_1(\chi)+\chi\act (e\tens w)=0,
\end{equation*}
and so $\TwIso$ maps $J^1_0(\E\tens\W)$ into $J^1_0(\E)\tens\W$.
\end{proof}
The operator $\TwIso$ formalises the process of twisting a first order
operator (on a $P$-module $\E$) by the twistor connection on $\W$. We apply
this to the exterior derivative in section~\ref{tdr}.

\section{Parabolic geometries}\label{pg}

Parabolic geometries can be described as Cartan geometries of type $(\g,P)$
where $\g$ is semisimple and the Lie algebra $\p$ of $P$ is a parabolic
subalgebra, i.e., a subalgebra containing a maximal solvable subalgebra of
$\g$. We need a few facts about parabolic subalgebras, all of which are
straightforward: we refer to~\cite{BE,CS,Slovak,Tanaka} for proofs.

The parabolic subalgebra splits naturally as the semidirect sum of a reductive
subalgebra $\g_0$ and a nilpotent ideal $\T\dual$, where $\T$ is the
orthogonal complement of $\p$ in $\g$ with respect to the Killing form of
$\g$---the nilpotent part of $\p$ is identified with $\T\dual$ using this
Killing form. Because of this duality, $\p\dual:=\g_0\subnormal\T$ is also a
parabolic subalgebra of $\g$. By choosing a Cartan subalgebra of the
semisimple part of $\g_0$ and extending it to a Cartan subalgebra of $\g$
inside $\g_0$, one can show (in the complexified setting) that parabolic
subalgebras of semisimple Lie algebras correspond, up to isomorphism, to
Dynkin diagrams with crossed nodes, where a node is crossed if the
corresponding root lies in the centre of $\g_0$. Real forms are classified in
a similar way (using, for instance, Satake diagrams). The distinction between
real and complex geometries does not cause any difficulties at this level:
some of the statements in the following require minor modification in the real
case, but we make little further comment on this.

We note that $[\g_0,\T]\subseteq \T$, $[\g_0,\T\dual]\subseteq \T\dual$ and
hence $\T$ and $\T\dual$ may be decomposed into graded nilpotent algebras by
the action of the centre of $\g_0$. This centre is nontrivial: in particular
there exists an element $E$ in the centre of $\g_0$ such that $\ad E$ has
positive integer eigenvalues on $\T$ and negative integer eigenvalues on
$\T\dual$, which may be normalized by the requirement that $1$ is an
eigenvalue. If $E$ acts by a scalar on a $\g_0$-module (as it does on an
irreducible $\g_0$-module), then this scalar will be called the
\emphdef{geometric weight}. By decomposing into irreducibles we can talk about
the geometric weights of any semisimple $\g_0$-module, and hence of any
element or function with values in that module. An important observation in
parabolic geometry is that although the grading of a $\g$ or $\p$-module by
geometric weight is not $\lie{p}$-equivariant, it induces a $\p$-equivariant
filtration. The associated graded vector space is the corresponding
$\g_0$-module.

If $P$ is a Lie group with Lie algebra $\p$ then we define $G_0$ to be the
subgroup $\{p\in P:\Ad_p(\g_0)\leq\g_0\}$; this has Lie algebra $\g_0$.  We
need to restrict the freedom in the choice of $P$ by assuming throughout that
$P=G_0\exp\T\dual$. This holds automatically if $G$ is a Lie group with Lie
algebra $\g$ and $P=\{g\in G:\Ad_g(\p)\leq\p\}$. The reason for this
assumption is that if we need to show a manifestly $G_0$-invariant
construction is $P$-invariant, we only need to check $\T\dual$-invariance. We
refer to such a $(\g,P)$, satisfying in addition the assumptions of the first
section, as a parabolic pair.

\begin{defn} A \emphdef{parabolic geometry} on $M$ is a Cartan geometry whose
type is a parabolic pair $(\g,P)$. If $\T$ is abelian, then this is called the
\emphdef{abelian} or \emphdef{almost Hermitian symmetric} case. A parabolic
geometry is said to be \emphdef{semiregular} if the geometric weights of the
curvature $\kappa$ are all nonpositive, and \emphdef{regular} if they are all
negative.
\end{defn}
In the abelian case, the centre is one dimensional, and the geometric weight
determines the action of the centre on an irreducible $\g_0$-module.  In fact
$\T$ itself has geometric weight $1$, and so an abelian parabolic geometric is
regular. Note that $\Lambda^2\T\dual\tens\p$ has negative geometric weights
(at most $-2$), so the (semi)regularity condition is a condition on the
torsion alone. Regularity ensures that the Lie bracket of vector fields on $M$
is compatible with the Lie bracket in $\T$---see~\cite{CS,Slovak,Tanaka}.

In practice, parabolic geometries are defined in terms of more primitive data,
which has to be prolonged (i.e., differentiated) to obtain the Cartan
geometry. It is natural to impose a further constraint on the curvature of
Cartan connections arising in this way, see~\ref{normal}. Here we give some
examples of geometric structures inducing such ``normal'' parabolic
geometries.

\subsection*{Conformal geometry}

It is well known that conformal geometry in $n\geq3$ dimensions (or M\"obius
geometry in dimension two~\cite{DMJC1}) can be described by a Cartan geometry
with $\g\isom\so(n+1,1)$. We fix a Lorentzian vector space $V$ of signature
$(n+1,1)$. Then the space of null lines in $V$ is the $n$-sphere viewed as a
conformal manifold, and the Lorentzian transformations act conformally. We
choose a point in $S^n$ and denote its tangent space, which is a conformal
vector space, by $\T$. The isotropy group fixing this null line is isomorphic
to $\CO(\T)\subnormal\T\dual$, with the conformal group $\CO(\T)$ acting on
$\T\dual$ in the obvious way. The Lorentzian Lie algebra $\g$, which is
semisimple for any $n\geq1$, is a vector space direct sum
$\g=\T\dsum\co(\T)\dsum\T\dual$. The geometric weight is the conformal weight.

Possible choices for $P$ are $\CO(\T)\subnormal\T\dual$, where $\CO(\T)$ may
or may not include the orientation reversing transformations, or
$\CSpin(\T)\subnormal\T\dual$. These parabolic geometries are called
(oriented) conformal geometry and conformal spin geometry respectively.

A more primitive definition of a conformal manifold is a manifold equipped
with an $L^2$-valued metric on the tangent bundle, where $L^1$ is the density
line bundle. We shall briefly describe how that Cartan connection arises
geometrically. A conformal manifold has a distinguished family of torsion-free
connections called~\emphdef{Weyl connections}, which form an affine space on
the space of $1$-forms. We can define the bundle of Weyl geometries $\cW$ as
the bundle of splittings of $J^1TM\to TM$ determined by the Weyl
connections. This is an affine bundle modelled on $T\dual M$ and the Weyl
connections are its sections. The principal bundle $\cG$ is the pullback of
$\cW$ to the bundle of conformal frames $\CO(M)$. The Cartan connection arises
from the observation that a $0$-jet of a section of $\cW$ can be extended
uniquely to a $1$-jet of a section with vanishing Ricci tensor.
Usually a more algebraic description is given: for a more detailed discussion,
with proofs, see~\cite{Baston1,CSS2,Ochiai}.

\medbreak

We describe the following examples even more briefly, our aim being only
to indicate the scope of parabolic geometry.

\subsection*{Projective geometry}
This is a parabolic geometry of type
$\bigl(\sgl(n+1,\R),\GL(n,\R)\subnormal\R^n\bigr)$. The structure is purely
second order, being given by a projective equivalence class of torsion-free
connections on the tangent bundle.

\subsection*{Quaternionic geometry}
This is a parabolic geometry in $n=4m$ dimensions of type
$\bigl(\sgl(m+1,\HQ),S(GL(1,\HQ)\times\GL(m,\HQ))\subnormal\HQ^m\bigr)$.  A
manifold is equipped with an (almost) quaternionic structure iff there is a
chosen rank $3$ Lie subalgebra bundle of $\End(TM)$ pointwise isomorphic to
the imaginary quaternions. A quaternionic structure is an almost quaternionic
structure with vanishing torsion.

\subsection*{Projective contact geometry}
There is a a contact geometry associated with each semisimple Lie
algebra. Projective contact geometry is a parabolic geometry of type
$\bigl(\symp(2(m+1),\R),\Symp(2m,\R)\subnormal\cH^{2m+1}\bigr)$, where
$\cH^{2m+1}$ is the Heisenberg group with Lie algebra $\R^{2m}\dsum\R$, the
Lie bracket being the standard symplectic form on $\R^{2m}$. A projective
contact manifold turns out to be a contact manifold together with a chosen
class of ``projectively equivalent'' partial connections.

\subsection*{CR geometry}
The geometry induced on strictly pseudoconvex hypersurfaces in complex
manifolds is a parabolic geometry of type
$\bigl(\su(m+1,1),\CU(m)\subnormal\cH^{2m+1}\bigr)$, where the Heisenberg Lie
algebra is now identified with $\C^m\dsum\R$ and the Lie bracket is the
imaginary part of the standard Hermitian form on $\C^n$. In fact a partial
integrability condition on an almost CR structure suffices to define the
Cartan geometry~\cite{CS}.

\subsection*{Quaternionic CR geometry}
We define quaternionic CR geometry to be a parabolic geometry of type $\bigl(
\symp(m+1,1),\R^+\times\Symp(1)\Symp(m)\subnormal\tilde\cH^{4m+3}\bigr)$,
where $\tilde\cH^{4m+3}$ is the Lie group of the nilpotent Lie algebra
structure on $\HQ^m\dsum\R^3$ whose Lie bracket is the direct sum of the three
symplectic forms on $\HQ^m$.

\subsection*{Pfaffian systems in five variables}
One of the first nontrivial Cartan geometries discovered (by Cartan, of
course~\cite{Cinq}) turns out to be an exceptional parabolic geometry.  One
starts from a Lie algebra $\T$ with basis $\{X_1,Y_1,Z_2,X_3,Y_3\}$ such that
$[X_1,Y_1]=Z_2$, $[X_1,Z_2]=X_3$, $[Y_1,Z_2]=Y_3$ and all other brackets are
trivial. Here the subscripts denote the geometric weight. A derivation of this
algebra is determined by its action on $X_1$ and $Y_1$, so the derivations
preserving geometric weight form a Lie algebra isomorphic to $\gl(2,\R)$ and
$E$ is the identity matrix. A more lengthy computation shows that there is
only a two dimensional space of derivations from $\T$ to
$\T\normalsub\gl(2,\R)$ lowering the geometric weight by one.  Prolonging
twice more gives a Lie algebra structure on $\g=\T\dsum\gl(2,\R)\dsum\T\dual$,
which turns out to be a real form of the exceptional Lie algebra $\g_2$.
Hence if one equips a $5$-manifold $M$ with a rank two subbundle $H$ of the
tangent bundle such that Lie brackets of vector fields in $H$ generate a rank
three subbundle, and that further Lie brackets with $H$ generate $TM$, then
one obtains, by~\cite{CS,Tanaka}, a parabolic geometry of type
$\bigl(\g_2,GL(2,\R)\subnormal\hat\cH^5\bigr)$ where $\hat\cH^5$ is the Lie
group of $\T$.

This example is in some sense typical: most ``normal'' parabolic geometries
are determined by a Pfaffian system on the tangent bundle~\cite{Yamaguchi}.
The preceding examples (apart from quaternionic CR geometry) are unusual in
this respect: the geometric structure involves additional data.

\smallbreak

In the main example of conformal geometry, we noted that Weyl geometries are
closely related to the Cartan principal bundle. Turning this around, we can
define a \emphdef{Weyl connection}, for an arbitrary parabolic geometry, to be
a $G_0$-equivariant section of $\cG\to\cG_0$ where $\cG_0$ is the principal
$G_0$-bundle $\cG/\exp\T\dual$. Since $\cG$ is a principal
$\exp\T\dual$-bundle over $\cG_0$, such a section is equivalently given by a
$P$-equivariant map $\cG\to\exp\T\dual$, i.e., a section of the associated
bundle $\cW:=\cG\cross_P\exp\T\dual\cong\cG/G_0$ over $M$, where $P$ acts on
$\exp\T\dual$ by the regular representation, not the adjoint
representation. Hence this is an affine bundle, the \emphdef{bundle of Weyl
geometries}, modelled on $T\dual M=\cG\cross_P\T\dual$---note that the adjoint
action of $P$ on $\T\dual$ is equivalent to its adjoint action on
$\exp\T\dual$. For any $P$-module $\E$, a Weyl connection (i.e., a section of
$\cW$ over $M$) identifies $\cG\cross_P\E$, filtered by geometric weight, with
$\cG_0\cross_{G_0}\E$, which is the associated graded bundle. We shall use
this observation later: for further information, see~\cite{CSl}.

\begin{defn}[\textsl{Parabolic twistors}]
Let $\W$ be a $(\g,P)$-module. Then $\T\dual$ acts on $\W$ and the image
$\T\dual\act\W$ of this action is a $P$-subrepresentation since $\T\dual$ is an
ideal of $\p$.  Hence there is a natural $P$-equivariant map
$\W\to\W_{\T\dual}$ where $\W_{\T\dual}:=\W/(\T\dual\act\W)$ is the space of
\emphdef{coinvariants} of $\T\dual$ acting on $\W$. Consequently, on a
parabolic geometry, sections $s$ of $W$ induce sections of $W_{T\dual
M}:=\cG\cross_P\W_{\T\dual}$. Parallel sections of $W$ will be called
\emph{parabolic twistors} and the induced sections of $W_{T\dual M}$ will be
called \emphdef{parabolic twistor fields}.
\end{defn}
Missing from this description is a \emphdef{twistor operator}: we shall see
later that there is a differential operator acting on sections of $W_{T\dual
M}$ which characterizes the parabolic twistor fields. The twistor operator is
the first operator in the curved BGG sequence which we shall construct. To do
this we need some Lie algebra homology.

\section{Lie algebra homology and cohomology}\label{lahc}

Any Lie algebra $\lie{l}$ possesses naturally defined homology and cohomology
theories with coefficients in an arbitrary $\lie{l}$-module $\W$. These
homology and cohomology groups can be constructed using a (projective or
injective) resolution of $\W$ by a Koszul complex of $\W$-valued alternating
forms. We shall apply this to parabolic geometries by taking $\W$ to be a
$\g$-module and letting $\lie{l}=\T$ or $\T^*$, using the vector space direct
sum $\g=\T\dsum\g_0\dsum\T\dual$. We focus on $\Lambda\T\dual\tens\W$, which
leads to the homology of $\T\dual$ or the cohomology of $\T$ with values in
$\W$.

\subsection*{$\T\dual$ homology with values in $\W$}

We first interpret the space $\Lambda^k\T\dual\tens\W$ as the space
$C_k(\T\dual,\W)$ of \emphdef{$k$-chains} on $\T\dual$ with values in $\W$.
This space carries a natural representation of $\p$: the action
on $\W$ is the restriction of the $\g$ action, while on $\Lambda^k\T\dual$ we
have:
\begin{equation}\label{dotdef}
\xi\act\beta := \sum_i \Lieb{\xi,\eps^i} \wedge (e_i\inner\beta)
\end{equation}
for $\xi\in\p$, where $e_i$ denotes a basis of $\T$ with dual basis
$\eps^i$. In the abelian case this action of $\T\dual\subseteq\p$ on
$\Lambda^k\T\dual$ is trivial. In general it is compatible with exterior
multiplication by $\alpha\in\T\dual$ in the sense that:
\begin{equation}\label{dotwedge}
\xi\act(\alpha\wedge\beta)=\alpha\wedge(\xi\act\beta)
+\Lieb{\xi,\alpha}\wedge\beta.
\end{equation}
There is also a compatibility with interior multiplication by $v\in\T$:
\begin{align}\label{dotinner}
v\inner(\xi\act\beta)
&=\xi\act(v\inner\beta)+\ip{\Lieb{\xi,\eps^i},v}e_i\inner\beta\\
\tag*{and so}
\xi\act(v\inner\beta)&=v\inner(\xi\act\beta)+\Lieb{\xi,v}\low_\T\inner\beta,
\end{align}
where $\Lieb{\xi,v}\low_{\T}$ denotes the $\T$ component of the Lie bracket in
$\g$, which is the coadjoint action of $\xi$ on $v\in\T\leq\g\dual$,
or equivalently the natural action on $v\in\T\isom\g/\p$.

Next we define the boundary operator or codifferential $\delTd$, where the
subscript denotes the Lie algebra which effectively acts in the following
definition:
\begin{equation}\begin{split}\label{codifdef}
\delTd\from C_k(\T\dual,\W)&\to C_{k-1}(\T\dual,\W)\\
\delTd( \beta \tens w )&=\sum_i\bigl(\tfrac12\eps^i\act(e_i\inner\beta)\tens w
+ e_i\inner\beta \tens \eps^i\act w\bigr).
\end{split}\end{equation}
For $k=0,1$ this definition means $\delTd w = 0$ and
$\delTd(\alpha \tens w) = \alpha\act w\,$.

\begin{lemma}[\textsl{Cartan's identity}] For $\alpha\in\T\dual$,
$\beta\in\Lambda^k\T\dual$, $w\in\W$ and $c\in C_k(\T\dual,\W)$ we have
\begin{align*}
\sum_ie_i\inner(\alpha\wedge\beta)\tens\eps^i\act w+
\sum_i\alpha\wedge(e_i\inner\beta)\tens\eps^i\act w&=\beta\tens\alpha\act w\\
\sum_i\eps^i\act\bigl(e_i\inner(\alpha\wedge\beta)\bigr)\tens w+
\sum_i\alpha\wedge\eps^i\act(e_i\inner\beta)\tens w
&=2\alpha\act\beta\tens w\qquad\qquad\\
\tag*{\textit{and consequently}} \delTd(\alpha\wedge c)+\alpha\wedge(\delTd c)
&=\alpha\act c\,.
\end{align*}
\proofof{lemma} The first part is immediate from the fact that
$e_i\inner(\alpha\wedge\beta)=\alpha(e_i)\beta-\alpha\wedge(e_i\inner\beta)$.
For the second part, we compute, using~\eqref{dotdef} and~\eqref{dotwedge},
\begin{align*}
\sum_i\eps^i\act(e_i\inner\alpha\wedge\beta)\tens w
&=\alpha\act\beta\tens w-\sum_i\eps^i\act(\alpha\wedge e_i\inner\beta)\tens w\\
&=\alpha\act\beta\tens w-\sum_i\bigl(\Lieb{\eps^i,\alpha}\wedge(e_i\inner\beta)
+\alpha\wedge\eps^i\act(e_i\inner\beta)\bigr)\tens w\displaybreak[0]\\
&=2\alpha\act\beta\tens w-\sum_i\alpha\wedge\eps^i\act(e_i\inner\beta)\tens w.
\end{align*}
The final formula (Cartan's identity) follows from the first two by
taking $c=\beta\tens w$.
\end{proof}
Cartan's identity perhaps best explains the curious factor of $1/2$ in
the definition of the codifferential. This factor is also crucial in
the following.
\begin{lemma}[\textsl{Boundary property}] The Lie algebra codifferential
defines a complex, i.e., $\delTd\circ\delTd=0$.
\proofof{lemma} If $\W$ is a trivial representation, then the definition of
the codifferential reduces to the term $\sum_i\tfrac12\eps^i\act
(e_i\inner\beta)\tens w$. We first show that the square of this term vanishes
separately by virtue of~\eqref{dotinner} and the Jacobi identity for
$\T\dual$:
\begin{align*}
\sum_{i,j} \eps^i\act(e_i\inner\eps^j\act(e_j\inner\beta))
&=\sum_{j,k}\Lieb{\eps^j,\eps^k}\act(e_k\inner e_j\inner\beta)
+\sum_{i,j}\eps^i\act\bigl(\eps^j\act(e_i\inner e_j\inner\beta)\bigr)\\
&=\frac12\sum_{i,j}\Lieb{\eps^i,\eps^j}\act(e_j\inner e_i\inner\beta)
=\frac12\sum_{i,j,k} \Lieb{\Lieb{\eps^i,\eps^j},\eps^k}\wedge
(e_k\inner e_j\inner e_i\inner\beta)\\
&=0.
\end{align*}
We can now compute $\delTd\circ\delTd$ in general:
\begin{align*}
\delTd\bigl(\delTd(\beta\tens w)\bigr)
&=\frac14\sum_{i,j} \eps^i\act(e_i\inner\eps^j\act(e_j\inner\beta))\tens w
+\sum_{i,j} e_i\inner e_j\inner\beta\tens\eps^i\act\eps^j\act w\\
&\quad+\frac12\sum_{i,j}e_i\inner\eps^j\act(e_j\inner\beta)\tens\eps^i\act w
+\frac12\sum_{i,j}\eps^i\act(e_i\inner e_j\inner\beta)\tens\eps^j\act w\\
&=0+\frac12\sum_{i,j} e_i\inner e_j\inner\beta\tens\Lieb{\eps^i,\eps^j}\act w\\
&\quad+\frac12\sum_{i,j}\bigl(e_i\inner\eps^j\act(e_j\inner\beta)
-\eps^j\act(e_i\inner e_j\inner\beta)\bigr)\tens\eps^i\act w,
\end{align*}
which vanishes by~\eqref{dotinner}.
\end{proof}
\begin{lemma}[\textsl{$\p$-equivariance}]\label{pact} For
$\alpha\in\T\dual$ and $c\in C_k(\T\dual,\W)$,
$\delTd(\alpha\act c)=\alpha\act (\delTd c)$.
\proofof{lemma} Both sides equal $\delTd\bigl(\alpha\wedge(\delTd
c)\bigr)$ by the previous two lemmas.
\end{proof}
It follows that $\delTd$ is $\p$-equivariant (since it is clearly
$\g_0$-equivariant).

\begin{defn}[\textsl{Homology}]
The cycles, boundaries and homology of $\delTd$ are given by:
\begin{align*}
Z_k(\T\dual,\W)
&:= \kernel \delTd\from C_k(\T\dual,\W)\to C_{k-1}(\T\dual,\W),\\
B_k(\T\dual,\W)
&:= \image \delTd\from C_{k+1}(\T\dual,\W)\to C_k(\T\dual,\W),\\
H_k(\T\dual,\W)
&:= Z_k(\T\dual,\W)/B_k(\T\dual,\W).
\end{align*}
\end{defn}
\noindent These are all $\p$-modules, and by Cartan's identity, $\T\dual$ maps
$Z_k(\T\dual,\W)$ into $B_k(\T\dual,\W)$ and hence acts trivially on the
homology $H_k(\T\dual,\W)$.

Note that the zero homology $H_0(\T\dual,\W)$ is the space of
coinvariants of $\W$ with respect to $\T\dual$, since in that case
$Z_0(\T\dual,\W)=\W$ and $B_0(\T\dual,\W)=\image(\act\from\T\dual\tens\W\to\W)
=\T\dual\act\W$, i.e., $H_0(\T\dual,\W)=\W/(\T\dual\act\W)=\W_{\T\dual}$.

\begin{examples}\label{confhom} The homology for the trivial representation
in the abelian case gives back the usual multilinear forms:
\begin{align*}
\W              &= \R,\\
H_k(\T\dual,\W) &= \Lambda^k\T\dual.
\end{align*}
In the case of conformal geometry, the Lorentzian vector space
$V=L^1\dsum\T^0\dsum L^{-1}$, where $\T^0=L^1\tens\T\dual$, is itself a
$\g$-module:
\begin{align*}
\W              &= L^1 \dsum \T^0 \dsum L^{-1},\\
H_0(\T\dual,\W) &= L^1,\\
H_k(\T\dual,\W) &= \Lambda^k\T\dual\cartan\T^0,\\
H_n(\T\dual,\W) &= \Lambda^n\T\dual\tens L^{-1}.
\end{align*}
Here elements in the \emphdef{Cartan product} $\Lambda^k\T\dual\cartan\T\dual$
can be thought of as tensors in $\Lambda^k\T\dual \tens \T\dual$ which are
alternating-free and trace-free (i.e., in the kernel of the natural maps to
$\Lambda^{k+1}\T\dual$ and $\Lambda^{k-1}\T\dual$)---this is the component
generated by the tensor product of highest weight vectors (in the complexified
representations if necessary).

Similarly for the adjoint representation $\W=\g$ we find:
\begin{align*}
\W &= \T\dsum \co(\T) \dsum\T\dual,\\
H_0(\T\dual,\W) &= \T,\\
H_1(\T\dual,\W) &= \T\dual \cartan \T,\\
H_k(\T\dual,\W) &= \Lambda^k \T\dual \cartan \so(\T),\\
H_{n-1}(\T\dual,\W) &= \Lambda^{n-1} \T\dual \cartan \T\dual,\\
H_n(\T\dual,\W) &= \Lambda^n \T\dual \tens \T\dual.
\end{align*}
Again elements in the Cartan product $\Lambda^k \T\dual \cartan \so(\T)$ can
be thought of as elements in the usual tensor product which are in the kernel
of natural $\co(\T)$-equivariant contractions and alternations. Elements in
the homology for $k=0,1,2$ have immediate geometric interpretations as
vectors, linearized conformal metrics and Weyl curvature tensors.
\end{examples}

\subsection*{$\T$ cohomology with values in $\W$}

Secondly we view $\Lambda^k\T\dual\tens\W$ as the space $C^k(\T,\W)$ of
\emphdef{$k$-cochains} on $\T$ with values in $\W$. From this point of view it
carries a natural $\p\dual$-action, where the action of $\chi\in\p\dual$ on
$\beta\in \Lambda^k\T\dual$ is
\begin{equation*}
\chi\act\beta= \sum_i \eps^i\wedge\bigl(\Lieb{e^i,\chi}\inner\beta\bigr)=
\sum_i \Lieb{\chi,\eps^j}\low_{\T\dual}\wedge(e_j\inner\beta),
\end{equation*}
with $\Lieb{\chi,\eps^j}\low_{\T\dual}$ denoting the $\T\dual$ component of the
Lie bracket.

The coboundary operator or differential $\dT\from C^k(\T,\W)\to
C^{k+1}(\T,\W)$ is defined by
\begin{equation*}
\dT( \beta\tens w )= \sum_i\bigl(\tfrac12\eps^i\wedge(e_i\act\beta)\tens w
+\eps^i\wedge\beta\tens e_i\act w\bigr).
\end{equation*}
This formula is minus the transpose of the formula for a boundary operator.
To be precise, it means that $\dT=-(\delT)\dual$, where $\delT$ is the
codifferential for $\T$ homology with values in $\W\dual$, whose $k$-chains
are $C_k(\T,\W\dual)=C^k(\T,\W)\dual$. It immediately follows that $\dT$
defines a complex and is $\p\dual$-equivariant (since $\delT$ is equivariant
with respect to $\p\dual=\g_0\dsum\T$ not $\p=\g_0\dsum\T\dual$). Cartan's
identity becomes $\dT(v\inner c)+v\inner\dT c=v\act c$ for $v\in\T$ and the
cohomology $H^k(\T,\W)$ of $\dT$ is naturally a $\p\dual$-module with $\T$
acting trivially.

\subsection*{Duality} In the above treatment of cohomology we made use
of the fact that it is dual to homology. More precisely, the $\T$ cohomology
with values in $\W$ is dual to the $\T$ homology with values in $\W\dual$:
\begin{align*}
C^k(\T,\W)\dual&=C_k(\T,\W\dual),&(\dT)\dual&=-\delT,&
H^k(\T,\W)\dual&=H_k(\T,\W\dual).\\
\intertext{Similarly, $\T\dual$ homology with values in $\W$ (the
first homology theory above) is dual to $\T\dual$ cohomology with
values in $\W\dual$:}
C_k(\T\dual,\W)\dual&=C^k(\T\dual,\W\dual),&(\delTd)\dual&=-\dTd,&
H_k(\T\dual,\W)\dual&=H^k(\T\dual,\W\dual).
\end{align*}
There is also a kind of Poincar\'e duality between $\T\dual$ homology and
cohomology (and similarly for $\T$): if $\T\dual$ is $n$-dimensional then
$C^k(\T\dual,\W\dual)\isom\Lambda^n\T\tens C_{n-k}(\T\dual,\W\dual)$ and
one easily checks that this intertwines boundary and coboundary so that
$H^k(\T\dual,\W\dual)\isom\Lambda^n\T\tens H_{n-k}(\T\dual,\W\dual)$.

We are interested primarily in $\delTd$ and $\dT$, both of which are
defined on $\Lambda\T\dual\tens\W$. It is natural to ask how they are
related. Since $\g$ is semisimple, one can use a Cartan involution to find
positive definite inner products on $\g$ and $\W$ with respect to
which $\delTd$ is minus the adjoint of $\dT$. From this, one obtains
Kostant's Hodge decomposition~\cite{Kostant}:
\begin{equation*}
\Lambda\T\dual\tens\W=\image\dT\dsum\bigl(\kernel\dT\intersect\kernel\delTd)
\dsum \image\delTd.
\end{equation*}
Furthermore, $\kernel\dT\intersect\kernel\delTd$ may be identified with the
kernel of Kostant's quabla operator $\Quabla=\delTd\dT+\dT\delTd$. The first
two terms in the Hodge decomposition give $\kernel\dT$ and the last two terms
give $\kernel\delTd$ and hence
\begin{equation*}
H^k(\T,\W)\isom\kernel\Quabla\isom H_k(\T\dual,\W).
\end{equation*}
This isomorphism is an isomorphism of $\g_0$-modules, although the cohomology
is viewed as a $\p\dual$-module with $\T$ acting trivially, whereas the
homology is viewed as a $\p$-module with $\T\dual$ acting
trivially. Similarly, the Hodge decomposition, as a direct sum, is only
$\g_0$-invariant, although the filtration
$0\leq\image\dT\leq\kernel\dT\leq\Lambda\T\dual\tens\W$ is $\p\dual$-invariant
and the filtration
$0\leq\image\delTd\leq\kernel\delTd\leq\Lambda\T\dual\tens\W$ is
$\p$-invariant.

\subsection*{The Main Theorem}

If $M$ is a parabolic geometry of type $(\g,P)$ and $\W$ is a $(\g,P)$-module,
then the Lie algebra homology groups are all $P$-modules and hence induce
vector bundles on $M$. We shall write $H_k(W)$ for $\cG\cross_P
H_k(\T\dual,\W)$, where $W=\cG\cross_P\W$. We also write $\Cinf(H_k(W))$ as a
shorthand for $\Cinf(M,H_k(W))$; more formally, we can interpret it as the
sheaf $U\mapsto\Cinf(U,H_k(W))$ for open subsets $U$ of $M$. Since $U$ is a
parabolic geometry in its own right, this slight of hand is more apparent than
real.

Our goal in the next two sections is to prove an explicit version of the
following result, the first part of which is due to \v Cap, Slov\'ak and
Sou\v cek~\cite{CSS4}, although our construction will be less complicated.

\begin{thm} Let $(M,\cform)$ be a parabolic geometry of type $(\g,P)$ and let
$\W$ be a finite dimensional $(\g,P)$-module. Then there is a naturally
defined sequence
\begin{equation*}
\Cinf(H_0(W))\;\stackrel{\OpInv_0}{\to}\;
\Cinf(H_1(W))\;\stackrel{\OpInv_1}{\to}\;
\Cinf(H_2(W))\;\stackrel{\OpInv_2}{\to}\;\cdots
\end{equation*}
of linear differential operators such that the kernel of the first operator is
isomorphic to the parabolic twistors associated to $W$ and the symbols of the
differential operators depend only on $(\g,P,\W)$ not $(M,\cform)$. If $M$ is
flat then this is sequence is locally exact and hence computes the cohomology
of $M$ with coefficients in the locally constant sheaf of parallel sections of
$W$.

Suppose further that $\W_1$, $\W_2$ and $\W_3$ are three finite dimensional
$(\g,P)$-modules with a nontrivial $(\g,P)$-equivariant linear map
$\W_1\tens\W_2\to\W_3$ (for instance $\W_3=\W_1\tens\W_2$). Then there are
nontrivial bilinear differential pairings
\begin{align*}
\Cinf(H_k(W_1)) &\cross \Cinf(H_\ell(W_2)) &\to&\quad
\Cinf(H_{k+\ell}(W_3))\\ (\alpha&,\beta)&\mapsto&\qquad \alpha\cpinv\beta
\end{align*}
whose symbols depend only on $(\g,P,\W_1,\W_2,\W_3)$ and which have the
following properties if $M$ is flat: for $k=\ell=0$ the pairing extends the
given pairing of twistors $\W_1\tens\W_2\to\W_3$, while more generally the
following Leibniz rule holds
\begin{equation*}
\OpInv_{k+\ell} (\alpha\cpinv\beta) 
= (\OpInv_k\alpha) \cpinv\beta + (-1)^k\alpha \cpinv (\OpInv_\ell\beta),
\end{equation*}
and hence the pairing descends to a cup product in cohomology.
\end{thm}

\section{The twisted deRham sequence}\label{tdr}

On a parabolic geometry $M$ of type $(\g,P)$ the chain complex of the previous
section induces, provided $\W$ is a $(\g,P)$-module, a complex of vector
bundles on $M$. If $W$ is the bundle induced by $\W$, then the bundle induced
by the $k$-chains is $\Lambda^kT\dual M\tens W$. The codifferential $\delTd$
induces a codifferential $\delTdM$ on the $k$-chain bundles, since it is
$P$-equivariant. On the other hand, $\dT$ is not $P$-equivariant and so does
not define an operator on the $k$-chain bundles on $M$. There is, however, a
first order differential operator, namely the exterior covariant derivative
(twisted deRham differential)
\begin{equation*}
\dtw\from J^1(\Lambda^kT\dual M\tens W)\to
\Lambda^{k+1}T\dual M\tens W
\end{equation*}
which behaves in many ways like $\dT$. To construct $\dtw$ formally, recall
that the invariant derivative defines an isomorphism from $J^1(\Lambda^kT\dual
M\tens W)$ to $\cG\cross_P J^1_0(\Lambda^k\T\dual\tens\W)$, which in turn
is isomorphic to $\cG\cross_P J^1_0(\Lambda^k\T\dual)\tens\W$. Hence we
need to find the $P$-equivariant map
$J^1_0(\Lambda^k\T\dual)\to\Lambda^{k+1}\T\dual$ induced by the exterior
derivative. So let $\alpha\from\cG\to\Lambda^k\T\dual$ be $P$-equivariant.
Then $\cform$ identifies $\alpha$ with a horizontal $P$-equivariant $k$-form on
$\cG$. Since exterior differentiation commutes with pullback, we can compute
the exterior derivative on $\cG$. This gives the following formula for
$d\alpha\from\cG\to\Lambda^{k+1}\T\dual$.
\begin{align*}
d\alpha(\xi_0,\ldots\xi_k)&=\sum_i(-1)^i(\partial_{\cform^{-1}(\xi_i)}\alpha)
(\xi_0,\ldots \hat\xi_i,\ldots\xi_k)\\ &\quad
+\sum_{i<j}(-1)^{i+j}\alpha(\cform[\cform^{-1}(\xi_i),\cform^{-1}(\xi_j)],
\xi_0,\ldots \hat\xi_i,\ldots\hat\xi_j,\ldots\xi_k)\\
\tag*{and so}
d\alpha&=\sum_i \eps^i\wedge\InvDer_{e_i}\alpha
-\sum_{i<j}\eps^i\wedge\eps^j\wedge
\bigl(\kappa(e_i,e_j)\inner\alpha\bigr)
+\tfrac12\sum_i\eps^i\wedge e_i\act\alpha\,.
\end{align*}
Combining this with the formula for $\TwIso$ in Proposition~\ref{ltt}, gives
the following result.

\begin{prop}[\textsl{Formal exterior derivatives}]\label{ecd} Let $\W$ a
$(\g,P)$-module. Then the exterior covariant derivative induces the
$P$-equivariant map
\begin{align}\notag
\qquad\sigma(\dtw)\from J^1_0(\Lambda^k\T\dual\tens\W)&
\to\Lambda^{k+1}\T\dual\tens\W\\ \label{extder}
(\phi_0,\phi_1)&\mapsto\sum_i\eps^i\wedge\phi_1(e_i)+\dT\phi_0-\sum_{i<j}
\eps^i\wedge\eps^j\wedge\bigl(\kappa(e_i,e_j)\inner\phi_0\bigr).\!\!\!\!\!
\intertext{The last term vanishes if the parabolic geometry is torsion-free.
In general it is $P$-equivariant, and so the $P$-equivariant map}\notag
\sigma(\dinv)\from J^1_0(\Lambda^k\T\dual\tens\W)&
\to\Lambda^{k+1}\T\dual\tens\W\\ \label{torsder} (\phi_0,\phi_1)&\mapsto
\sum_i\eps^i\wedge\phi_1(e_i)+\dT\phi_0
\end{align}
induces a exterior covariant derivative $\dinv$ with torsion \textup(unless
$\cform$ is torsion-free\textup).
\end{prop}
Thus we see that although the zero order operator $\dT$ is not
$P$-equivariant, we can correct it by a first order term: the symbol of the
exterior derivative (wedge product). There are two ways to do this. The
simplest formula~\eqref{torsder} (with no torsion correction) gives an
exterior covariant derivative with torsion, but an extra term can be added to
remove the torsion~\eqref{extder}.  On the bundle level, these exterior
derivatives are related by
\begin{equation*}
d^\g s=d^\cform s-\sum_{i<j}\eps^i\wedge\eps^j\wedge K\low_M(e_i,e_j)\inner s
\end{equation*}
for any section $s$ of $\Lambda^kT\dual M\tens W$: note that only the torsion
part of $K\low_M$ contributes to the contraction. This difference is also
illustrated by the following result.
\begin{prop}[\textsl{Curvature}] The composites $\Rtw=(\dtw)^2$
and $\Rinv:=(\dinv)^2$ acting on a section $s$ of $\Lambda^kT\dual
M\tens W$ are given by
\begin{align*}
\Rtw s&=\trace \bigl(X\mapsto K\low_M\wedge X\act s\bigr)\\
\Rinv s&=\trace\bigl(X\mapsto -K\low_M\wedge\InvDer_X s\bigr).
\end{align*}
Here $X\in\gM$: in the first formula, the $\gM$-values of
$K\low_M$ act on the $W$-values of $s$, while in the second formula, the
$\gM$-values are contracted with the invariant derivative.
\proofof{prop} The first formula is clear: the square of the $d^\g$ is the
wedge product action of the curvature of $\TwConn$. For the second
formula, let $f\from\cG\to\Lambda^k\T\dual\tens\W$ be $P$-equivariant. Then
\begin{equation*}
\dinv f=\sigma(\dinv)(f,\InvDer f)=\dT f
+\sum_i\eps^i\wedge\InvDer_{e_i}f.
\end{equation*}
We must apply $\dinv$ to this expression. To do this, note that $\dT$
is constant on $\cG$, so that
\begin{equation*}
\sum_j\eps^j\wedge\InvDer_{e_j} \dinv f
=\sum_j\eps^j\wedge\dT\InvDer_{e_j}f
+\sum_{i,j}\eps^j\wedge\eps^i\wedge\InvDer_{e_j}(\InvDer_{e_i}f)
\end{equation*}
and therefore
\begin{align*}
(\dinv)^2f&=(\dT)^2 f+\sum_i\dT(\eps^i\wedge\InvDer_{e_i}f)+
\sum_j\eps^j\wedge\dT\InvDer_{e_j}f
+\sum_{i,j}\eps^j\wedge\eps^i\wedge\InvDer_{e_j}(\InvDer_{e_i}f)\\
&=0+\frac12\sum_{i,j}\bigl(\eps^j\wedge\Lieb{e_j,\eps^i}\low_{\T\dual}
\wedge\InvDer_{e_i} f
+\eps^j\wedge\eps^i\wedge
(\InvDer_{e_j}(\InvDer_{e_i}f)-\InvDer_{e_i}
(\InvDer_{e_j}f))\bigr)\displaybreak[0]\\
&=\frac12\sum_{i,j}\eps^j\wedge\eps^i\wedge\bigl(
\InvDer_{e_j}(\InvDer_{e_i}f)
-\InvDer_{e_i}(\InvDer_{e_j}f)
-\InvDer_{\Lieb{e_j,e_i}}f\bigr)\\
&=\frac12\sum_{i,j}\eps^j\wedge\eps^i\wedge\InvDer_{\kappa(e_i,e_j)}f
=-\sum_i\ip{\zeta^i,\kappa}\wedge\InvDer_{\chi_i}f,
\end{align*}
where $\chi_i$ is a basis of $\g$ with dual basis $\zeta^i$.
\end{proof}
These vanish if $K$ is zero, or if $K$ takes values in a subspace of
$\p$ acting trivially on $\W$.

\section{The BGG sequence and cup product}\label{bgg}

The key tool for proving the main theorem is a family of differential
operators
\begin{equation*}
\PiInv_k\from\Cinf(\Lambda^kT\dual M\tens W)\to
\Cinf(\Lambda^kT\dual M\tens W)
\end{equation*}
which vanish on $\image\delTdM$, map into $\kernel\delTdM$, and induce the
identity on homology. As motivation for the construction of such an operator,
recall Kostant's quabla operator $\Quabla=\delTd\dT+\dT\delTd$ (with
$\kernel\Quabla\isom H_k(\T\dual,\W)$) and the Hodge decomposition:
\begin{equation*}
\Lambda^k\T\dual\tens\W=\image\dT\dsum\kernel\Quabla\dsum\image\delTd\,.
\end{equation*}
The projection onto $\kernel\Quabla$ in this direct sum has image
contained in $\kernel\delTd$ and induces the identity on
homology. Unfortunately, it is not $\p$-equivariant. Ignoring this for the
moment, note that $\Quabla$ is invertible on its own image and so the
projection onto $\kernel\Quabla$ may be written $\iden-\Quabla^{-1}\Quabla$. A
more refined formula may be obtained by observing that $\Quabla$ commutes with
$\dT$, and hence so does $\Quabla^{-1}$ on the image of $\Quabla$. Therefore:
\begin{equation*}
\iden-\Quabla^{-1}(\delTd\dT+\dT\delTd)
=\iden-\Quabla^{-1}\delTd\dT-\dT\Quabla^{-1}\delTd\,.
\end{equation*}
The advantage of this formula is that we only need the inverse of
$\Quabla$ on $\image\delTd$ (which is a $\p$-module). Indeed, we only
need the operator $\Quabla^{-1}\delTd$.

We now address the problem of $\p$-equivariance. Of course $\Quabla$ fails to
be $\p$-equivariant simply because $\dT$ is not $\p$-equivariant.  However, in
the previous section we noted that one resolution is to replace $\dT$ with a
first order differential operator: either $\dinv$ or $\dtw$. We shall
concentrate first on the former, but return to the latter at the end of
the section.

\begin{defn}[\textsl{First order quabla operator}] Let $M$ be a parabolic
geometry of type $(\g,P)$ and let $\W$ be a $(\g,P)$-module. Then the
\emph{quabla operator} on $\Lambda T\dual M\tens W$ is the first order
differential operator $\QuaInv=\delTdM\dinv+\dinv\delTdM$.
\end{defn}

Note that $\QuaInv$ commutes with $\delTdM$ and also maps $k$-forms to
$k$-forms, so it preserves $B_k(W)=\image\delTdM\from\Lambda^{k+1} T\dual
M\tens W\to\Lambda^k T\dual M\tens W$. In the flat case it also commutes with
$\dinv$, but in general $\QuaInv\circ\dinv-\dinv\circ\QuaInv=
\delTdM\circ\Rinv-\Rinv\circ\delTdM$.

\begin{thm}\label{inv} Suppose $M$ is a parabolic geometry of type $(\g,P)$
and $\W$ is a finite dimensional $(\g,P)$-module.  Then
$\QuaInv\from\Cinf(B_k(W))\to\Cinf(B_k(W))$ is invertible. Furthermore
the inverse is a differential operator of finite order.

\proofof{thm} To prove that $\QuaInv$ has a two-sided differential inverse, we
choose a Weyl connection, i.e., a section of the bundle of Weyl geometries
$\cW$. Such a section always exists locally on $M$---in the smooth, rather
than analytic, category, they exist globally, since $\cW$ is an affine bundle,
but we shall only need local sections, since we are constructing a local
operator and, by the uniqueness of two-sided inverses, the local inverses
patch together. Hence we assume we have a section over all of $M$, which
identifies the $k$-chain bundles, filtered by geometric weight, with the
associated graded bundles. The operators $\Quabla$ and $\dT$, which are
$G_0$-invariant, but not $P$-invariant, define operators on associated graded
bundles, and hence, using the Weyl connection, on the $k$-chain bundles
themselves.
\begin{lemma} $\Quabla^{-1}(\QuaInv-\Quabla)
\from\Cinf(B_k(W))\to\Cinf(B_k(W))$ is nilpotent.
\proofof{lemma}[Proof of Lemma] Since $\W$ is finite dimensional, the
$\p$-module $B_k(\T\dual,\W)$ decomposes into finitely many irreducible
$\g_0$-submodules and clearly the action of $\T\dual$ lowers the geometric
weight. Suppose that $s\from M\to B_k(W)$ takes values in an subbundle
associated to an  irreducible $\g_0$-submodule. Now
\begin{align*}
(\QuaInv-\Quabla)s&=(\delTdM(\dinv-\dT)+(\dinv-\dT)\delTdM)s\\
&=\sum_i\bigl(\delTdM(\eps^i\wedge\InvDer_{e_i}s)
+\eps^i\wedge\delTdM\InvDer_{e_i}s\bigr)\\
&=\sum_i\eps^i\act\InvDer_{e_i}s
\end{align*}
which has lower geometric weight, since the covariant derivative preserves the
filtration, while the action of $T\dual M$ lowers the weight.  Finally,
$\Quabla^{-1}$ is $\g_0$-equivariant and so it preserves the geometric
weight. Hence $\Quabla^{-1}(\QuaInv-\Quabla)$ lowers the geometric weight, so
it is nilpotent.
\end{proof}
Writing $\cN=-\Quabla^{-1}(\QuaInv-\Quabla)$, we have
$\QuaInv=\Quabla(\iden-\cN)$. Therefore the two-sided inverse
$\QuaInv^{-1}$ is given by a Neumann series:
\begin{equation*}
\QuaInv^{-1}=(\iden-\cN)^{-1}\Quabla^{-1}
=\biggl(\sum_{k\geq0}\cN^k\biggr)\Quabla^{-1}.
\end{equation*}
This inverse is a differential operator whose order is the degree of
nilpotency of the first order differential operator $\cN$. Hence it is an
inverse on any open subset of $M$. Our construction involved the choice of
Weyl connection, but the inverse constructed is of course independent of this
choice.
\end{proof}
We can now define differential operators $\Qinv$ and $\PiInv$
from $\Cinf(\Lambda T\dual M\tens W)$ to itself:
\begin{equation*}
\Qinv=\QuaInv^{-1}\delTdM\,,\qquad \PiInv=\iden-\dinv\circ
\Qinv-\Qinv\circ\dinv.
\end{equation*}
Clearly $\PiInv$ preserves the degree of the form. This gives our operators
$\PiInv_k$.

\begin{rem}\label{same} The first order quabla operator $\QuaInv$
maps sections of $Z_k(W)$ into $B_k(W)$. This means in particular that it
preserves $B_k(W)$, but also that it descends to an operator on
$C_k(W)/Z_k(W)$.  In the construction of $\Qinv$ and $\PiInv$, we only used
the invertibility of $\QuaInv$ on $B_k(W)$, but in fact it is also invertible
on $C_k(W)/Z_k(W)$: in the algebraic setting this holds for Kostant's
$\Quabla$ and exactly the same Neumann series argument goes through. Hence one
might prefer to define $\tilde\Qinv=\delTdM\QuaInv^{-1}$, where $\QuaInv^{-1}$
is now the inverse on $C_k(W)/Z_k(W)$ and $\tilde\Qinv$ acts on $C_k(W)$ by
first passing to the quotient---clearly this is well defined since
$Z_k(W)=C_k(W)\intersect\ker\delTdM$ by definition.

Now observe that $\Qinv$ and $\tilde\Qinv$ are both given by isomorphisms from
$C_k(W)/Z_k(W)$ to $B_{k-1}(W)$. Composing on each side by $\QuaInv$ gives
$\QuaInv\Qinv\QuaInv=\delTdM\QuaInv$ and
$\QuaInv\tilde\Qinv\QuaInv=\QuaInv\delTdM$. Since $\QuaInv$ is an
isomorphism on $C_k(W)/Z_k(W)$ and on $B_{k-1}(W)$, and it commutes with
$\delTdM$, we deduce that $\tilde\Qinv=\Qinv$.
\end{rem}

We now establish the fundamental properties of $\PiInv$.

\begin{prop}[\textsl{Calculus of $\Pi$-operators}]\label{calc} The operator
$\PiInv_k\from\Cinf(\Lambda^kT\dual M\tens W)\to\Cinf(\Lambda^kT\dual M\tens
W)$ has the following properties.
\begin{enumerate}
\item $\PiInv_k$ vanishes on $\image\delTdM\,$\textup:\quad
$\PiInv_k\circ\delTdM=0$.
\item $\PiInv_k$ maps into $\kernel\delTdM\,$\textup:\quad
$\delTdM\circ\PiInv_k=0$.
\item On $\kernel\delTdM$, $\PiInv_k\cong\iden\mod
\image\delTdM$, i.e., $\PiInv_k$ induces the identity on homology.
\item $\dinv\circ\PiInv_k-\PiInv_{k+1}\circ \dinv=\Qinv\circ\Rinv-\Rinv\circ
\Qinv$.
\item $(\PiInv_k)^2=\PiInv_k+\Qinv\circ \Rinv\circ \Qinv$ and so $\PiInv_k$
is a projection in the flat case, and for $k=0$.
\item $\PiInv\circ\QuaInv=\Qinv\circ\Rinv\circ\delTdM$ and
$\QuaInv\circ\PiInv=\delTdM\circ\Rinv\circ \Qinv$.
\end{enumerate}
Thus in the flat case $\PiInv_k$ is a differential projection onto a subspace
of $\kernel\delTdM$ complementary to $\image\delTdM$ and is a chain map from
the deRham complex to itself; $\Qinv$ is a chain homotopy between $\PiInv_k$
and $\iden$.
\proofof{prop} The first three results follow from $\kernel\delTdM=\kernel
\Qinv$ and $\image \Qinv=\image\delTdM$ (since $\QuaInv^{-1}$ is the inverse
on $\image\delTdM$). The fourth fact follows easily from the definition of
$\PiInv$ and using this, the fifth fact is an immediate calculation:
\begin{align*}
(\PiInv_k)^2&=\PiInv_k(\iden-\dinv\circ \Qinv)=\PiInv_k
-(\dinv\circ\PiInv_{k-1}-\Qinv\circ \Rinv
+\Rinv\circ \Qinv)\circ \Qinv\\
&=\PiInv_k+\Qinv\circ \Rinv\circ\Qinv.
\end{align*}
The last part also follows easily from the definition of $\PiInv$.
\end{proof}
The first two properties allow us to define two further operators:
\begin{align}
\PiInv_k\circ\represent &\from \Cinf(H_k(W))
\to\Cinf(\Lambda^kT\dual M\tens W)\label{crepr}\\
{\project}\circ\PiInv_k&\from\Cinf(\Lambda^kT\dual M\tens W)
\to\Cinf(H_k(W))\label{cproj}
\end{align}
where $\project$ denotes the projection from the kernel of $\delTdM$ to
homology and $\represent$ means the choice of a representative of the homology
class. Thus $\PiInv_k\circ\represent$ gives a canonical differential
representative for homology classes and ${\project}\circ\PiInv_k$ is a
canonical differential projection onto homology. By Proposition~\ref{calc}
(vi), the canonical representative is the unique representative in
$\kernel\QuaInv$, whereas the canonical projection vanishes on
$\image\QuaInv$.

The first operator~\eqref{crepr} was originally constructed by Baston in the
abelian case~\cite{Baston2}, and \v Cap, Slov\'ak and Sou\v cek in
general~\cite{CSS4}. Note that on $\kernel\delTdM$,
$\PiInv=\iden-\QuaInv^{-1}\delTdM\dinv$. It is interesting to note that in the
flat abelian case, Baston obtained the Neumann series formula for this
in~\cite{Baston2}, equation (8).

We define $\OpInv_k={\project}\circ\PiInv_{k+1}\circ
\dinv\circ\PiInv_k\circ\represent$. Since $\dinv\circ\PiInv_k$ maps
$\kernel\delTdM$ into $\kernel\delTdM$ already, this equals ${\project}\circ
\dinv\circ\PiInv_k\circ\represent$, so we only actually need~\eqref{crepr}.
For the pairings we really need~\eqref{crepr} and~\eqref{cproj} to define
\begin{equation*}
\cpinv={\project}\circ\PiInv_{k+\ell}\circ\wedge
\circ(\PiInv_k\circ\represent,\PiInv_\ell\circ\represent)
\end{equation*}
where $\wedge$ denotes wedge product of forms contracted by the
pairing $W_1\tens W_2\to W_3$.

The main theorem is now straightforward (apart from the independence result
for the symbols---see the appendix): in the flat case we have a locally exact
resolution because $\PiInv$, as a chain map on the deRham resolution by
sheaves of smooth sections, is homotopic to the identity, and the Leibniz rule
follows from the corresponding Leibniz rule for the wedge product. In the
curved case we have the following results.

\begin{prop}[\textsl{Composition}]\label{compcurv}
$\OpInv_{k+1}\circ\OpInv_k={\project}\circ\PiInv_{k+2}\circ\Rinv\circ
\PiInv_k\circ\represent$.
\proofof{prop} By definition
$\OpInv_{k+1}\circ\OpInv_k={\project}\circ\dinv\circ\PiInv_{k+1}\circ
\dinv\circ\PiInv_k\circ\represent$. Now commute $\dinv_{k+1}$ past
$\PiInv_{k+1}$ using the $\Pi$-operator calculus of the previous proposition.
\end{proof}

\begin{prop}[\textsl{Leibniz rule}]\label{leib} For $\alpha\in\Cinf(H_k(W_1))$
and $\beta\in\Cinf(H_\ell(W_2))$,
\begin{multline*}
\OpInv_{k+\ell}(\alpha\cpinv \beta)=
\OpInv_k\alpha\cpinv \beta +(-1)^k \alpha\cpinv \OpInv_\ell \beta\\
+\Bigl[\PiInv_{k+\ell+1}\Bigl(\bigl(\Qinv \Rinv \PiInv_k \alpha\bigr)
\wedge \PiInv_\ell \beta +(-1)^k \PiInv_k \alpha \wedge
\bigl(\Qinv \Rinv \PiInv_\ell \beta\bigr)
-\Rinv \Qinv\bigl(\PiInv_k \alpha\wedge\PiInv_\ell \beta\bigr)\Bigr)\Bigr].
\end{multline*}
Here, and henceforth, we write $[\ldots]$ for the projection
to homology, and $\PiInv_k$ for $\PiInv_k\circ\represent$.
\proofof{prop} This again follows easily from Proposition~\ref{calc}:
\begin{align*}
\OpInv_{k+\ell}(\alpha\cpinv \beta)
&=[\PiInv_{k+\ell+1}\dinv\PiInv_{k+\ell}
(\PiInv_k\alpha\wedge\PiInv_\ell\beta)]\\
&=[\PiInv_{k+\ell+1}\dinv(\PiInv_k \alpha\wedge\PiInv_\ell \beta)]
-[\PiInv_{k+\ell+1}\Rinv\Qinv(\PiInv_k \alpha\wedge\PiInv_\ell \beta)].
\end{align*}
The first term can be expanded using the Leibniz rule for the exterior
derivative:
\begin{equation*}
\dinv(\PiInv_k \alpha\wedge\PiInv_\ell \beta)=
\dinv\PiInv_k \alpha\wedge\PiInv_\ell \beta
+(-1)^k\PiInv_k \alpha\wedge\dinv\PiInv_\ell \beta.
\end{equation*}
We insert the projections $\PiInv_{k+1}$, $\PiInv_{\ell+1}$ using the
definition $\iden=\PiInv+\dinv\circ\Qinv+\Qinv\circ\dinv$. The first
correction term does not contribute, since $\dinv\PiInv_k \alpha$ and
$\dinv\PiInv_\ell \beta$ are in $\kernel\delTdM$, while the second correction
gives two further curvature terms as stated.
\end{proof}

We now consider the other choice of exterior covariant derivative: the
torsion-free operator $\dtw$. This change makes no difference if the parabolic
geometry is torsion-free. In the presence of torsion, we can construct what we
believe is a more natural curved analogue of the BGG complex, although to do
this, we need to assume that the parabolic geometry is regular, i.e., the
geometric weights of the curvature $\kappa$ are negative (which is a condition
on the torsion).  Under this assumption the extra torsion correction in the
formula for $\dtw$ does not cause any problems in the proof of nilpotency in
Theorem~\ref{inv}, since its action still lowers the geometric weight, and all
other details of the proofs are unchanged. We thus obtain operators $\QuaTw$,
$\Qtw$, $\PiTw_k$, $\OpTw_k$, $\cptw$ satisfying the same formulae with
$\dinv$ replaced by $\dtw$ and $\Rinv$ by $\Rtw$.

The ``torsion-free'' BGG sequences seem to us to be more natural, because
$\Rtw$ is always zero order, given by wedge product with the curvature
$K_M$. Another reason for preferring $\dtw$ is the differential Bianchi
identity.
\begin{prop} Let $M$ be a Cartan geometry of type $(\g,P)$ with curvature
$K\low_M\in\Cinf(\Lambda^2T\dual M\tens\gM)$. Then $\dtw K\low_M=0$.
\end{prop}
We combine $\dtw K\low_M=0$ with a well-known definition.
\begin{defn}\label{normal} A parabolic geometry is said to be \emphdef{normal}
if $\delTdM K\low_M=0$.
\end{defn}
\begin{thm} Let $(\cG,\cform)$ be a normal regular parabolic geometry of
type $(\g,P)$ on $M$. Then the curvature $K\low_M$ is uniquely determined
by its homology class $[K\low_M]$, via the formula
\begin{equation*}
K\low_M=\PiTw_2[K\low_M]
\end{equation*}
and $[K\low_M]$ therefore satisfies $\OpTw_2[K\low_M]=0$, where
$\OpTw_2$ is the operator in the torsion-free BGG sequence
associated to $\g$.

Furthermore the composite of two operators in the torsion-free BGG sequence
associated to $\W$ is given by
\begin{equation*}
\OpTw_{k+1}\OpTw_k\alpha=[K\low_M]\cptw\alpha,
\end{equation*}
where the cup product is contracted by the pairing $\g\tens\W\to\W$ given by
the $\g$-action.
\proofof{thm} $\PiTw_2[K\low_M]$ is the unique element of
$\kernel\delTdM\intersect\kernel\QuaTw$ whose homology class is
$[K\low_M]$. But $K\low_M$ itself satisfies $\dtw K\low_M=0$ and $\delTdM
K\low_M=0$ and hence $\QuaTw K\low_M=0$.

For the second part, observe that
\begin{equation*}
[K\low_M] \cptw \alpha
=\bigl[\PiTw_{k+2}\bigl(\PiTw_2[K\low_M]\wedge\PiTw_k \alpha\bigr)\bigr]
=\bigl[\PiTw_{k+2}(K\low_M\wedge\PiTw_k \alpha)\bigr]
=\bigl[\PiTw_{k+2}\Rtw\PiTw_k \alpha\bigr]
\end{equation*}
which is the composite ${\project}\circ\PiTw_{k+2}\circ\Rtw\circ\PiTw_k
\circ\represent$ acting on $\alpha$.
\end{proof}
For conformal geometry in four dimensions or more, the curvature of the Cartan
connection is obtained by applying a first order operator to the Weyl
curvature, as is well known. Even in the general context, the observation that
the curvature is uniquely determined by its (co)homology class is an old one:
see~\cite{CS,Tanaka}. Our approach reveals that the proofs in these references
appear technical because they amount to the construction of
$\PiTw_2\circ\represent$ in this special case. Also, by working with Lie
algebra homology, rather than cohomology, the explicit differential operator
reconstructing the full curvature is realized as an operator on $M$, rather
than $\cG$, as in the conformal case.

In the second part of this theorem, it is slightly awkward that the action of
$\g$ needs to be specified. There is a convenient device to make this happen
automatically. Suppose that all $(\g,P)$-modules of interest belong to the
symmetric algebra or the tensor algebra of $\W=\W_1\dsum\W_2\dsum\cdots$; for
instance, $\W$ could be the standard representation or the direct sum of the
fundamental representations of $\g$. Now work either in the universal
enveloping algebra of $\W\idealsub\g$ (with trivial bracket on $\W$), or in
the tensor algebra of $\W\dsum\g$ modulo the ideal generated by
$X\tens\phi-\phi\tens X-X\act\phi$ for $X\in\g$ and $\phi\in\W\dsum\g$ (a
semiholonomic enveloping algebra---see the appendix). This algebra is filtered
by finite dimensional $\g$-modules, where the action is induced by the action
on $\g\dsum\W$, \textit{and definitely \emph{not} by left multiplication with
elements of $\g$}.  It follows that for any differential form $\alpha$ with
values in the associated algebra bundle,
$\Rtw\alpha=K\low_M\wedge\alpha-\alpha\wedge K_M$, where the curvature $K_M$
is viewed as a $2$-form with values in the copy of $\gM$ in this algebra
bundle. The properties of $\cK=[K_M]$ and $\OpTw$ established in the above
theorem may now be rewritten:
\begin{equation}\label{BianCurv}
\OpTw_2 \cK=0,\qquad\qquad
\OpTw_{k+1}\OpTw_k\alpha=\cK\cptw\alpha-\alpha\cptw\cK.
\end{equation}
The curvature terms in the Leibniz rule may be rewritten in a similar way:
\begin{equation}\label{Leib}
\OpTw_{k+\ell}(\alpha\cptw\beta)=\OpTw_k\alpha\cptw \beta +(-1)^k
\alpha\cptw\OpTw_\ell \beta-\Ip{\cK,\alpha,\beta}
+\Ip{\alpha,\cK,\beta}-\Ip{\alpha,\beta,\cK},
\end{equation}
where the triple products are defined by
\begin{equation*}\notag
\Ip{\cK,\alpha,\beta}=\Bigl[\PiTw_{k+\ell+1}\Bigl(
\PiTw_2\cK\wedge\Qtw\bigl(\PiTw_k\alpha\wedge\PiTw_\ell \beta\bigr)
-\Qtw\bigl(\PiTw_2\cK\wedge\PiTw_k\alpha\bigr)\wedge\PiTw_\ell\beta
\Bigr)\Bigr]
\end{equation*}
and similarly for the other two products, although the first term acquires a
sign $(-1)^k$. The contractions with $\cK$ happen automatically in this
combination of triple products. If $\alpha$ and $\beta$ belong to BGG
subsequences associated to $\W_1,\W_2$ then the formula can be contracted
further using any $(\g,P)$-equivariant linear map $\W_1\tens\W_2\to\W_3$.

These triple products may seem ad hoc, but in fact this is the first
appearance of natural trilinear differential pairings closely related to
Massey products. For any $(\g,P)$-equivariant linear map
$\W_1\tens\W_2\tens\W_3\to\W_4$, one can define a trilinear differential
pairing from $\Cinf(H_k(W_1))\times\Cinf(H_\ell(W_2))\times\Cinf(H_m(W_3))$ to
$\Cinf(H_{k+\ell+m-1}(W_4))$ by
\begin{equation*}
\Ip{\alpha,\beta,\gam}=\Bigl[\PiTw_{k+\ell+m-1}\Bigl(
(-1)^k\PiTw_k\alpha\wedge\Qtw\bigl(\PiTw_\ell \beta\wedge\PiTw_m \gam\bigr)
-\Qtw\bigl(\PiTw_k \alpha\wedge\PiTw_\ell \beta\bigr)\wedge\PiTw_m \gam
\Bigr)\Bigr].
\end{equation*}
This measures the failure of the cup product to be associative: one
may compute that
\begin{align}\label{trip}
\OpTw_{k+\ell+m-1}\Ip{\alpha,\beta,\gam}
&=(\alpha\cptw \beta)\cptw \gam-\alpha\cptw(\beta\cptw \gam)\\
&\quad-\Ip{\OpTw_k \alpha,\beta,\gam}-(-1)^k\Ip{\alpha,\OpTw_\ell \beta,\gam}
-(-1)^{k+\ell}\Ip{\alpha,\beta,\OpTw_m \gam}\notag\\
&\quad+\Ip{\cK,\alpha,\beta,\gam}-\Ip{\alpha,\cK,\beta,\gam}
+\Ip{\alpha,\beta,\cK,\gam}-\Ip{\alpha,\beta,\gam,\cK}\notag
\end{align}
where the quadruple products each have five terms. In the flat case, this
formula verifies that the cup product is associative in BGG cohomology. Note,
though, that in practice, one often destroys this associativity by using
incompatible (nonassociative) pairings to define the cup products: it is
crucial above that the same map $\W_1\tens\W_2\tens\W_3\to\W_4$ is used for
$(\alpha\cptw \beta)\cptw\gam$ and $\alpha\cptw(\beta\cptw \gam)$.

The relation, in the flat case, with a Massey product is as follows: if
$\OpTw_k\alpha=\OpTw_\ell \beta=\OpTw_m\gam=0$ and if also $\alpha\cptw
\beta=\OpTw_{k+\ell-1}A$ and $\beta\cptw\gam=\OpTw_{\ell+m-1}C$, then
\begin{equation*}
\OpTw_{k+\ell+m-1}\bigl(A\cptw\gam-(-1)^k\alpha\cptw C
-\Ip{\alpha,\beta,\gam}\bigr)=0
\end{equation*}
and hence we obtain a partially defined triple product of BGG cohomology
classes, with an ambiguity coming from the choice of $A$ and $C$.  Again the
role of $\Ip{\alpha,\beta,\gam}$ is to correct the failure of $\cptw$ to be
associative on the BGG cochain complex. Note that the two terms in
$\Ip{\alpha,\beta,\gam}$ modify the lifts $\PiTw A$ and $\PiTw C$ of $A$ and
$C$ from Lie algebra homology.

\section{Curved $A_\infty$-algebras}\label{Ainf}

The formulae~\eqref{BianCurv},~\eqref{Leib} and~\eqref{trip} give the first
four defining relations of a curved $A_\infty$-algebra. In the case of
vanishing curvature, such algebras were introduced by Stasheff~\cite{Stash}
nearly forty years ago. An $A_\infty$-algebra is a graded vector space $A$
equipped with a sequence of multilinear maps $\mu_k\colon{\tens}^kA\to A$ of
degree $2-k$ satisfying some identities. (In fact only parity really matters,
and $\mu_k$ has parity $k$ mod $2$.) In the original formulation, $\mu_0=0$,
$\mu_1$ is a differential, and $\mu_2$ is ``strongly homotopy associative'',
i.e., it is associative up to a homotopy given by $\mu_3$, which in turn
satisfies higher order associativity conditions. In the presence of curvature,
we require that for each $m\geq 0$,
\begin{equation*}\notag
\sum_{\substack{j+k=m+1\\j\geq1,k\geq0}}\sum_{\ell=0}^{j-1}
(-1)^{k+\ell+k\ell+k|\alpha_1\ldots\alpha_\ell|}
\mu_j\bigl(\alpha_1,\ldots\alpha_\ell,
\mu_k(\alpha_{\ell+1},\ldots\alpha_{k+\ell}),
\alpha_{k+\ell+1},\ldots\alpha_m\bigr)=0,
\end{equation*}
for all $\alpha_1,\ldots\alpha_m\in A$ of homogeneous degree, where
$|\alpha_1\ldots\alpha_\ell|$ denotes the sum of the degrees.  The usual
definition, with the sign conventions of~\cite{Markl}, is recovered by putting
$\mu_0=0$. For the Lie analogue of $L_\infty$-algebras, the general curved
case has been introduced by Zwiebach~\cite{Zwiebach} within the context of
String Field Theory, where the presence of $\mu_0$ is interpreted as a
non-conformal background, related to (genus $0$) vacuum vertices.  In our
setting, $\mu_0$ is the (background) curvature, and we now indicate briefly
how such a curved $A_\infty$-algebra structure arises.
Following~\cite{GM,JL,Merk}, we first work on the level of the chain bundles
and define $\lam_m$ inductively, for $m\geq 2$, by
\begin{equation*} 
\lam_m(a_1,\ldots a_m)=\sum_{\substack{j+k=m\\ j,k\geq 1}}
(-1)^{(k-1){\textstyle(}j+|a_1\ldots a_j|{\textstyle)}}
\Qtw\lam_j(a_1,\ldots a_j)\wedge\Qtw\lam_k(a_{j+1},\ldots a_m)
\end{equation*}
where we formally set $\Qtw\lam_1=-\iden$. Note that the number of terms in
$\lam_m$ is given by the Catalan number $\frac1{m+1}\binom{2m}{m}$.  On the
homology bundles, we then define: $\mu_0=\cK$, $\mu_1(\alpha_1)=\OpTw\alpha_1$
and $\mu_m(\alpha_1,\ldots\alpha_m)=
[\PiTw\lambda_m(\PiTw\alpha_1,\ldots\PiTw\alpha_m)]$ for $m\geq2$.

In order to prove that this is a curved $A_\infty$-algebra, it is convenient
to make use of the observation that an $A_\infty$-algebra structure on a
vector space $A$ is equivalently an odd coderivation of square zero on the
tensor coalgebra of $A$ (with the grading of $A$ shifted to get the signs
right)---see for instance~\cite{Markl}, Example 1.9. Although the tensor
coalgebra of the sheaf of sections of the homology bundles of an enveloping
algebra makes us a bit dizzy, we are only using this formalism as a way to
compute identities for multilinear differential operators which avoids dealing
with huge expressions and complicated signs.

To obtain the coderivation, put $\mu=\sum_{m\geq 0}\mu_m\colon\Tens A\to A$,
let $p_i\colon\Tens A\to{\tens}^i A$ be the projection, and define $\mu^c$ by
$p_0\mu^c=0$, $p_1\mu^c=\mu$ and $\Delta\circ\mu^c
=(\iden\tens\mu^c+\mu^c\tens\iden)\circ\Delta$, where
$\Delta(a_1\tens\ldots\tens a_k)=\sum_j(a_1\tens\dots\tens
a_j)\tens(a_{j+1}\tens\dots\tens a_k)$ is the coproduct. The defining
relations of an $A_\infty$-algebra are now equivalent to $(\mu^c)^2=0$,
although it suffices to check that $\mu\mu^c=0$, since $(\mu^c)^2$ is the
coderivation $(\mu\mu^c)^c$.

In our case, we have $\mu=\project\PiTw(K_M+\dtw+\lambda)\PiTw\represent$,
where $\lambda=\sum_{m\geq2}\lambda_m$ and $\PiTw$ is extended to the tensor
coalgebra as $\sum\PiTw\tens\cdots\tens\PiTw$. The proof that the induced
coderivation has square zero follows~\cite{Merk}, except that we must deal
with curvature terms. Such terms appear in five ways: the curvature explicitly
in the definition of $\mu$; the term $(\dtw)^2=\Rtw$ in $\mu\mu^c$; from
$\PiTw^2=\PiTw+\Qtw\Rtw\Qtw$; from $\project\PiTw\dtw\PiTw=
\project\PiTw(\dtw-\Rtw\Qtw)$; and from $\PiTw\dtw\PiTw\represent=
(\dtw-\Qtw\Rtw)\PiTw\represent$. Note that the shift in the grading changes
some signs and that the recursive definition of $\lambda_m$ is equivalent to
$\lambda= \lambda_2\bigl((\Qtw\lambda-p_1)\tens(\Qtw\lambda-p_1)\bigr)\Delta$.
Omitting the lift and projection to homology, we have
\begin{align*}
\mu\mu^c&=
\PiTw(K_M+\dtw+\lambda)\PiTw\bigl(\PiTw(K_M+\dtw+\lambda)\PiTw\bigr)^c\\&
=\PiTw\dtw\PiTw^2(K_M+\dtw+\lambda)\PiTw
+\PiTw\lambda\bigl(\PiTw^2(K_M+\dtw+\lambda)\bigr)^c\PiTw\\&
=\PiTw\Rtw p_1\PiTw+\PiTw(\dtw-\Rtw\Qtw)\lambda\PiTw
+\PiTw\lambda
\bigl(K_M+\dtw-\Qtw\Rtw p_1+\PiTw\lambda+\Qtw\Rtw\Qtw\lambda\bigr)^c\PiTw\\&
=\PiTw\bigl(\dtw\lambda+\lambda(\dtw)^c+\lambda(\PiTw\lambda)^c\bigr)\PiTw
+\PiTw\lambda_2
\bigl(K_M\tens(\Qtw\lambda-p_1)+(\Qtw\lambda-p_1)\tens K_M\bigr)\PiTw\\
&\qquad\qquad\qquad\qquad\qquad
+\PiTw\lambda_2\bigl((\Qtw\lambda-p_1)\tens(\Qtw\lambda-p_1)\bigr)
(K_M^c\tens\iden+\iden\tens K_M^c)\Delta\PiTw\\
&\qquad\qquad\qquad\qquad\qquad
+\PiTw\lambda\bigl(\Qtw\lambda_2\bigl(
K_M\tens(p_1-\Qtw\lambda)+(p_1-\Qtw\lambda)\tens K_M\bigr)\bigr)^c\PiTw\\
&=\PiTw\bigl(\lambda\lambda^c
+\dtw\lambda+\lambda(\dtw)^c-\lambda([\dtw,\Qtw]\lambda)^c\bigr)\PiTw\\
&\quad+\PiTw\lambda_2\bigl(\Qtw\lambda K_M^c\tens(\Qtw\lambda-p_1)
+(\Qtw\lambda-p_1)\tens \Qtw\lambda K_M^c\bigr)\Delta\PiTw\\
&\quad-\PiTw\lambda\bigl(\Qtw\lambda_2
\bigl(K_M\tens(\Qtw\lambda-p_1)+(\Qtw\lambda-p_1)\tens K_M\bigr)\bigr)^c\PiTw.
\end{align*}
Next we compute that $\lambda\lambda^c=
\lambda_2\bigl(\Qtw\lambda\lambda^c\tens(\Qtw\lambda-p_1)
+(\Qtw\lambda-p_1)\tens \Qtw\lambda\lambda^c\bigr)\Delta$---the term
$\lambda_2\bigl(
\lambda\tens(\Qtw\lambda-p_1)+(\Qtw\lambda-p_1)\tens\lambda\bigr)\Delta$
vanishes by expanding $\lambda$ and using the associativity of $\lambda_2$.
It follows by induction that $\lambda\lambda^c=0$. Similarly,
$\dtw\lambda+\lambda(\dtw)^c-\lambda([\dtw,\Qtw]\lambda)^c=\lambda_2\Bigl(
\Qtw\bigl(\dtw\lambda+\lambda(\dtw)^c-\lambda([\dtw,\Qtw]\lambda)^c\bigr)
\tens(\Qtw\lambda-p_1)+(\Qtw\lambda-p_1)\tens
\Qtw\bigl(\dtw\lambda+\lambda(\dtw)^c-\lambda([\dtw,\Qtw]\lambda)^c\bigr)
\Bigr)\Delta$, and so it also follows that
$\dtw\lambda+\lambda(\dtw)^c-\lambda([\dtw,\Qtw]\lambda)^c=0$.
One more recursive argument shows that the curvature terms cancel as well.

\begin{rem} J. Stasheff has pointed out to us that this sort of
result can also be obtained using the techniques of Homological Perturbation
Theory, at least in the flat case.  The crucial idea is that $\Qtw$ defines
strong deformation retraction data for the coderivation determined by
$\dtw$. The methods of~\cite{GLS} may then be used to transfer the
perturbation of this coderivation induced by the wedge product to the Lie
algebra homology bundles.
\end{rem}

Finally, we remark that restricting the above to the (super)symmetric
coalgebra gives an $L_\infty$-algebra, in which one can work with
$\W\idealsub\g$ instead of its enveloping algebra.

\section{The dual BGG sequences}\label{dualBGG}

The BGG cochain sequence of Lie algebra homology bundles $H_k(W)$ is dual to a
chain sequence of Lie algebra cohomology bundles, generalizing the deRham
chain complex of exterior divergences. To fix notations, recall that the
latter is a complex
\begin{equation*}
\Cinf(L^{-n})\;\stackrel{\delta}{\leftarrow}\;
\Cinf(L^{-n}\tens TM)\;\stackrel{\delta}{\leftarrow}\;
\Cinf(L^{-n}\tens TM)\;\stackrel{\delta}{\leftarrow}\;\cdots
\end{equation*}
where $L^{-n}$ is the oriented line bundle of densities and $\delta$ is the
exterior divergence, i.e., on vector field densities $\delta=\divg$, the
natural divergence, and in general it is adjoint to $d$ in the sense that for
$\alpha\in\Cinf(\Lambda^kT\dual M)$ and
$a\in\Cinf(L^{-n}\tens\Lambda^{k+1}TM)$, we have
\begin{equation*}
\divg(\alpha\capinner a)=\ip{d\alpha, a}+\ip{\alpha,\delta a},
\end{equation*}
where $\ip{\theta,\alpha\capinner a}=\ip{\theta\wedge\alpha,a}$ for any
$1$-form $\theta$. For compactly supported sections, the complex can be
augmented by $\int\colon\Cinf_0(L^{-n})\to\R$, giving a homology theory.

A simple way to obtain a dual BGG sequence is to twist the BGG sequence of the
dual $(\g,P)$-module $\W\dual$ by the flat line bundle of pseudoscalars
$L^{-n}\tens\Lambda^n TM$, where $n=\dim M$: this is the orientation line
bundle of $M$ and is associated to a one dimensional $(\g,P)$-module, on which
only $G_0\leq P$ might act nontrivially. By Poincar\'e duality for Lie algebra
(co)homology, such a twist amounts to replacing $H_k(W\dual)$ with
$L^{-n}\tens H^{n-k}(W\dual)$, where $H^{n-k}(W\dual)=H_{n-k}(W)\dual$.
Writing $\Opinv^{n-1-k}$ for this twist of $\OpInv_k$, we obtain a
sequence
\begin{equation*}\notag
\Cinf(L^{-n}\tens H^0(W\dual))\;\stackrel{\Opinv^0}{\leftarrow}\;
\Cinf(L^{-n}\tens H^1(W\dual))\;\stackrel{\Opinv^1}{\leftarrow}\;
\Cinf(L^{-n}\tens H^2(W\dual))\;\stackrel{\Opinv^2}{\leftarrow}\cdots
\end{equation*}
of linear differential operators. In the flat case, this is, by construction,
an injective resolution of the sheaf of parallel sections of $L^{-n}\tens
\Lambda^n TM\tens W\dual$, beginning with $\Opinv^{n-1}$, but it is
natural to view it instead as a projective resolution of the dual of the sheaf
of parallel sections of $W$ by working with compactly supported sections and
defining, for $a\in \Cinf_0(L^{-n}\tens H^0(W\dual))$, $\Ip{\int a,w}
=\int\ip{a,[w]}$ for any parallel section $w$ of $W$.

This point of view is further amplified by constructing the dual BGG sequences
directly from the sequence of exterior divergences twisted by the twistor
connection on $W\dual$. In the presence of torsion, there are two
possibilities: $\delinv=-(\dinv)^*$ or $\deltw=-(\dtw)^*$. As in the previous
section, we define $\hQuaInv=\dTd\delinv+\delinv\dTd$ and find that
it is invertible on $\Cinf(L^{-n}\tens B^k(W\dual))$ and $\Cinf(L^{-n}\tens
C^k(W\dual)/Z^k(W\dual))$, where $C^k(W\dual)=\cG\cross_P C^k(\T\dual,\W\dual)
=\Lambda^kTM\tens W\dual$. Hence we obtain operators $\hat\Qinv$
and $\hat\PiInv$. Furthermore, this construction is adjoint
to the construction of the previous section: since $\delinv=-(\dinv)^*$
and $\dTd=-(\delTd)^*$, we have $\hQuaInv=(\QuaInv)^*$,
$B^k(W\dual)=(C_k(W)/Z_k(W))\dual$, $C^k(W\dual)/Z^k(W\dual)=B_k(W)\dual$,
and hence, by Remark~\ref{same}, $\hat\Qinv=-(\Qinv)^*$, so that
$\hat\PiInv=(\PiInv)^*$.

The dual BGG operators obtained above by Poincar\'e duality are therefore
equivalently defined by
$\Opinv^k={\project}\circ\delinv\circ\hat\PiInv_k\circ\represent$.
Associated to a pairing $\W_1\tens\W_2\to\W_3$, the analogue of the cup
product is a ``cap product'' between cochains and chains:
\begin{align*}
\Cinf(H_k(W_1))\times&\Cinf(L^{-n}\tens H^{k+\ell}(W_3\dual))
&\to&\quad \Cinf(L^{-n}H^{\ell}(W_2\dual))\\
(\alpha,&b)&\mapsto&\qquad\quad\alpha\capinv b
\end{align*}
satisfying a Leibniz rule up to curvature terms. This can be defined
by twisting the cup product by $L^{-n}\tens\Lambda^nTM$ and using
Poincar\'e duality, or by the formula
\begin{equation*}
\capinv={\project}\circ\hat\PiInv_{k+\ell}\circ\capinner
\circ(\PiInv_k\circ\represent,\hat\PiInv_\ell\circ\represent)
\end{equation*}
where $\capinner$ denotes the contraction of forms with multivectors together
with the pairing $W_1\tens W_3\dual\to W_2\dual$. Here $\alpha\capinner a$ for
$\alpha\in\Cinf(\Lambda^kT\dual M)$ and
$a\in\Cinf(L^{-n}\tens\Lambda^{k+\ell}TM)$ is defined by
$\ip{\theta,\alpha\capinner a}=\ip{\theta\wedge\alpha,a}$ for any $\ell$-form
$\theta$.

The Leibniz rule, for $\alpha\in \Cinf(H_k(W_1))$ and
$b\in\Cinf(L^{-n}\tens H^{k+\ell}(W_3\dual))$, is:
\begin{multline*}
\Opinv^{\ell}(\alpha\capinv b) = \alpha\capinv(\Opinv^{k+\ell}
b)-(-1)^\ell(\OpInv_k\alpha)\capinv b\\
+\Bigl[\PiInv_{\ell-1}\Bigl(\PiInv_k\alpha \capinner
\bigl(\hat\Qinv\hat\Rinv\hat\PiInv_{k+\ell} b\bigr)
-(-1)^\ell\bigl(\Qinv\Rinv\PiInv_k\alpha\bigr)\capinner\hat\PiInv_{k+\ell} b
-\hat\Rinv\hat\Qinv\bigl(\PiInv_k\alpha\capinner\hat\PiInv_\ell
b\bigr)\Bigr)\Bigr].
\end{multline*}
Similar results hold for the torsion-free sequence
$\OpTw^k={\project}\circ\deltw\circ\hat\PiTw_k\circ\represent$ (in the
regular case) and if the parabolic geometry is also normal, the composite of
dual BGG operators is given by cap product with $[K_M]$ and the correction
terms to the Leibniz rule are given by triple products of $[K_M]$, $\alpha$
and $b$.

The cap product gives a neat way to see the duality between $\OpTw_k$ and
$\OpTw^k$. Consider the pairing $\W\tens\R\to\W$, with cap product
\begin{equation*}
\Cinf(M,H_k(W))\times\Cinf(M,L^{-n}\tens H^{k+\ell}(W\dual))\to
\Cinf(M,L^{-n}\tens H^{\ell}(\R)).
\end{equation*}
For $\ell=0$, $H^0(\R)=\R$ and this pairing is the duality pairing of $H_k(W)$
and $H^k(W\dual)$, tensored with the density bundle $L^{-n}$. When $\ell=1$,
$H^1(\R)$ is the subbundle of $TM$ associated to $\g_1$, the geometric weight
$1$ subspace of $\g$. Hence we obtain a pairing with values in vector
densities. The claim is that the Leibniz rule for $\ell=1$ becomes
\begin{equation*}
\divg\alpha\captw b=\ip{\OpTw_k\alpha,b}+\ip{\alpha,\OpTw^k b}
\end{equation*}
with no curvature corrections. The inner curvature corrections vanish because
$Z_k(W)=B^k(W\dual)^0$ and $Z^k(W\dual)=B_k(W)^0$, so that the contractions of
$\PiTw$ with $\hat\Qtw$ and $\hat\PiTw$ with $\Qtw$ are zero. The outer
correction vanishes because the curvature is acting on the trivial
representation. A similar result holds for the analogue of the triple product
rule~\eqref{trip}, showing that for adjoint pairings of $(\g,P)$-modules,
$c\mapsto\beta\captw c$ is adjoint to $\alpha\mapsto\alpha\cptw\beta$.

\section{General applications}\label{genapp}

We turn now to potential applications of the BGG sequence and cup product.  We
restrict ourselves first to general discussions: more explicit examples are
given in the next section.

\subsection*{Twistor operators} The exterior covariant derivatives
$\dtw$ and $\dinv$ on $W$-valued $0$-forms are both simply the covariant
derivative $\TwConn$ on sections of $W$. Also $\kernel\delTdM=W$ and so a
parabolic twistor $f$ is the natural representative $\PiTw_0[f]$ of its
homology class, which is a parabolic twistor field.  Hence
$\PiTw_0\circ\represent$ is a ``jet operator'' which assigns a parabolic
twistor to the corresponding parabolic twistor field. In the flat case, the
\emphdef{twistor operator}
$\OpTw_0=\project\circ\TwConn\circ\PiTw_0\circ\represent$ characterizes
parabolic twistor fields as solutions of a differential equation.  In the
curved case it is natural to define parabolic twistor fields by the kernel of
this operator, but $\OpTw_0\phi=0$ only implies that
$\TwConn\PiTw_0\phi=\Qtw\Rtw\PiTw_0\phi$, and so $\PiTw_0\phi$ might not be
parallel in general.

\subsection*{Twistor algebra} $\g$-modules form an algebra under direct
sum and tensor product. The cup product $\Cinf(M,H_0(W_1))
\times\Cinf(M,H_0(W_2))\to\Cinf(M,H_0(W_1\tens W_2))$ defines an algebra
structure on sections of the corresponding homology bundles. The Leibniz rule
shows that in the flat case the cup product algebra extends the algebra of
twistors. A similar observation can be made for $\g$-modules under Cartan
product (provided one is careful with identifications between isomorphic
representations)---in this case the cup product is zero order, given by the
Cartan product of zeroth homologies.

\subsection*{Deformation theory and moduli spaces} Suppose that $M$ is a
compact manifold admitting a flat parabolic geometry of type $(\g,P)$.  What
is the moduli space of flat parabolic geometries of type $(\g,P)$ on $M$?  A
first approximation to this question is to study deformations of the given
flat structure $(\cG,\cform)$.  We discuss briefly deformations of regular
normal parabolic geometries, with emphasis on deformations of flat structures.

Fixing the principal $P$-bundle $\cG\to M$, Cartan connections of type
$(\g,P)$ form an open subset of an affine space modelled on the
$P$-equivariant horizontal $1$-forms $T\cG\to\g$ (it is an open subset because
of the condition that $\cform$ is an isomorphism on each tangent
space). Therefore a small deformation of $\cform$ may be written
$\cform_\eps=\cform+\tilde\alpha_\eps$ where $\tilde\alpha_\eps$ is a curve of
such $P$-equivariant horizontal $1$-forms with $\tilde\alpha_0=0$. The
curvature of $\cform_\eps$ is $K^\eps(U,V)=d\cform_\eps(U,V)
+[\cform_\eps(U),\cform_\eps(V)]$ and passing to associated bundles gives
\begin{equation}\label{defcurv}
K_M^\eps=K_M+\dtw\alpha_\eps+\alpha_\eps\wedge\alpha_\eps,
\end{equation}
where $\alpha_\eps\colon TM\to\gM$ (with $\alpha_0=0$) and we think of $\gM$
as a subbundle of its universal enveloping algebra bundle, so that
$(\alpha_\eps\wedge\alpha_\eps)(X,Y)=
\alpha_\eps(X)\cdot\alpha_\eps(Y)-\alpha_\eps(Y)\cdot\alpha_\eps(X)=
[\alpha_\eps(X),\alpha_\eps(Y)]\in\gM$ for $X,Y\in TM$ (equivalently, we can
work directly with the Lie bracket in $\gM$). Differentiating with respect to
$\eps$ at $\eps=0$, gives, up to third order,
\begin{equation*}
\dot K\low_M=\dtw\dot\alpha,\qquad
\ddot K\low_M=\dtw\ddot\alpha+2\dot\alpha\wedge\dot\alpha,\qquad
\dddot K\low_M=\dtw\dddot\alpha+3(\ddot\alpha\wedge\dot\alpha+
\dot\alpha\wedge\ddot\alpha).
\end{equation*}

Suppose now that $K_M=0$. Then, in order for $\cform_\eps$ to be normal, we
need $\delTdM\dot K\low_M=0$. Also $\dtw\dot K\low_M=(\dtw)^2\dot\alpha=0$,
and so $\dot K\low_M=\PiTw_2[\dot K\low_M]$ and $\OpTw_2[\dot K\low_M]=0$.
Adding $\dtw s$ to $\dot\alpha$ does not alter $\dot K\low_M$, and so one can
assume $\delTdM\dot\alpha=0$. Hence $\QuaTw\dot\alpha=0$ and $\dot\alpha$
represents a homology class $A=[\dot\alpha]$. We then have
$\dot\alpha=\PiTw_1A$, $[\dot K\low_M] =[\dtw\dot\alpha]=\OpTw_1A$. This does
not completely fix the freedom to add $\dtw s$ to $\dot\alpha$: we can still
add $\OpTw_0 f$ to $A$.

To summarize, we see that the linearized theory is controlled by the BGG
complex with $\W=\g$: an infinitesimal deformation of $\cform$ (as a regular
normal parabolic geometry) is given by a section $A$ of $H_1(\gM)$ and
$\OpTw_1A$ is the linearized curvature. Since $\OpTw_1\OpTw_0f=0$,
$A=\OpTw_0f$ as just an infinitesimal gauge transformation. Hence the formal
tangent space to the moduli space is the first cohomology of the complex
\begin{equation*}\notag
\Cinf(M,H_0(\gM))\xrightarrow{\OpTw_0}\Cinf(M,H_1(\gM))\xrightarrow{\OpTw_1}
\Cinf(M,H_2(\gM))\xrightarrow{\OpTw_2}\Cinf(M,H_3(\gM)).
\end{equation*}
This is only the actual tangent space if all infinitesimal deformations can be
integrated. We first consider second order deformations.  If $\dot K\low_M=0$
(i.e., $\OpTw_1A=0$), then normality implies that $\delTdM\ddot K\low_M=0$,
and since also $\dtw\ddot K\low_M=0$, $\ddot K\low_M=\PiTw_2[\ddot
K\low_M]$. Hence it suffices to consider
$[\PiTw_2\dtw\ddot\alpha]+[\PiTw_2(\dot\alpha\wedge\dot\alpha)]$ and the
second term is $A\cptw A$.  As before, we can assume $\delTdM\ddot\alpha=0$,
so that the first term is $\OpTw_1[\ddot\alpha]$.  Hence we have a quadratic
obstruction to solving $\ddot K_M=0$: we need $A\cptw A$ to be in the image of
$\OpTw_1$. The Leibniz rule gives $\OpTw_2(A\cptw A)=2(\OpTw_1 A)\cptw A=0$
and so the obstruction is the class of $A\cptw A$ in the second cohomology of
above complex.

The obstructions to building a formal power series all lie in this second
cohomology space, but the construction involves the $A_\infty$-algebra of
multilinear operators, not just the cup product (alternatively we can work in
the $L_\infty$-algebra associated to $\g$). The reason for this can be seen at
third order: $\QuaTw\ddot\alpha=\delTdM\dtw\ddot\alpha$ is not zero in
general, and $\ddot\alpha=\PiTw_1\ddot\alpha
+\Qtw\dtw\ddot\alpha=\PiTw_1\ddot\alpha-2\Qtw(\dot\alpha\wedge\dot\alpha)$.
Hence if $A_1=[\dot\alpha]$ and $A_2=\frac12[\PiTw_1\ddot\alpha]$, we have
\begin{align*}\notag
\dot\alpha=\PiTw_1 A_1\qquad&\mathrm{and}\qquad\ddot\alpha=2\bigl(\PiTw_1A_2
-\Qtw(\PiTw_1A_1\wedge\PiTw_1A_1)\bigr),\\ \tag*{and therefore}
[\PiTw_2(\ddot\alpha\wedge\dot\alpha+\dot\alpha\wedge\ddot\alpha)]&=
2\bigl(A_2\cptw A_1+A_1\cptw A_2+\ip{A_1,A_1,A_1}\bigr).
\end{align*}
The Leibniz and triple product rules~\eqref{Leib},~\eqref{trip} imply that
this is in the kernel of $\OpTw_2$, and its cohomology class is the second
obstruction.

In general if $A=\sum_{j=1}^k A_j\eps^j$ satisfies the equation $\OpTw_1
A+A\cptw A+\ip{A,A,A}+\cdots=0$ to order $k$ in $\eps$, then $\OpTw_2(A\cptw
A+\ip{A,A,A}+\cdots)$ vanishes to order $k+1$ in $\eps$ and the cohomology
class of the degree $k+1$ term of $A\cptw A++\ip{A,A,A}+\cdots$ is the
obstruction to finding $A_{k+1}$ such that $\tilde A=\sum_{j=1}^{k+1}
A_j\eps^j$ satisfies the equation to order $k+1$ in $\eps$.

This deformation theory parallels numerous examples in algebra and geometry
which have been studied since~\cite{Ger,KS}. It would be interesting to extend
it to half-flat geometries such as selfdual conformal $4$-manifolds~\cite{DF}.

\subsection*{Linear field theories} In the flat case (when the BGG sequence
is a complex) it is natural to view it as a linear gauge theory. The beginning
of the sequence gives the kinematics, while the end gives the dynamics;
equivalently, the dynamics are given by the beginning of the dual BGG
sequence.
\begin{equation*}\notag
\begin{array}{cccccc}
\mathrm{charges}&\mathrm{gauges}&&\mathrm{potentials}&&
\mathrm{kinematic\;fields}\\
0\to\W\to\!\!\!&
 \Cinf(H_0(W))&\!\!\xrightarrow[\mathrmsl{Twistor}]{\OpTw_0}\!\!&
 \Cinf(H_1(W))&\!\!\xrightarrow[\mathrmsl{Potential}]{\OpTw_1}\!\!&
 \Cinf(H_2(W))\\[4mm]
0\leftarrow\!\W\dual\!\leftarrow\!\!\!&
 \Cinf(L^{-n}H^0(W\dual))&\!\!\xleftarrow[\mathrmsl{ConsLaw}]{\OpTw^0}\!\!\!&
 \Cinf(L^{-n}H^1(W\dual))&\!\!\xleftarrow[\mathrmsl{FieldEqn}]{\OpTw^1}\!\!\!&
 \Cinf(L^{-n}H^2(W\dual))\\
\!\!\mathrm{dual\; charges}\!\!\!\!\!&\mathrm{fluxes}&&\mathrm{sources}&&
\mathrm{dynamic\;fields}
\end{array}\end{equation*}
The prototype for such a sequence is the deRham complex describing
electromagnetism.  We shall present some further justification for this point
of view on linear field theory in the following, but we refer to~\cite{Thesis}
for a more thorough discussion. The assumption that the BGG sequence is a
complex means that potentials coming from a gauge via $\OpTw_0$ have no
kinematic field, and that sources coming from a dynamic field automatically
satisfy the conservation law $\OpTw^0$. Not shown is the kinematic field
equation $\OpTw_2$ which is automatically satisfied by kinematic fields coming
from a potential.

Given a $(\g,P)$-equivariant pairing $\W_1\tens\W_2\to\W_3$, we have in
particular a cup product
$\Cinf(M,H_0(W_1))\times\Cinf(M,H_k(W_2))\to\Cinf(M,H_k(W_3))$. This means
that twistors in $W_1$ can be used to transform objects in the field theory
associated to $\W_2$ to the field theory associated to $\W_3$.  The Leibniz
rule implies that this transformation will be compatible with the operators in
the sequence. This is often called the translation principle.

\subsection*{Relations to deRham complexes}

We restrict attention here to the abelian case, so the BGG sequence of the
trivial representation $\R$ is the deRham complex of exterior derivatives, and
the dual BGG sequence is the dual deRham complex of exterior divergences.
For any $\g$-module $\W$, we always have a pairing $\W\tens\R\to\W$, and so,
given a twistor in $W$, we can construct potentials and kinematic fields in
the $\W$-theory from differential forms using the cup product (``kinematic
helicity raising''), while the cap product gives multivector densities from
dynamic fields and sources (``dynamic helicity lowering''). Similarly, the
pairing $\W\dual\tens\R\to\W\dual$ shows that a twistor in $W\dual$ can be
used to construct dynamic fields and sources from multivector densities via
cap product (``dynamic helicity raising''), and differential forms from
potentials and kinematic fields via cup product (``kinematic helicity
lowering''). Let us focus on the cap product
\begin{equation*}
\Cinf(H_k(W))\times\Cinf(L^{-n}\tens H^{k+\ell}(W\dual))\to
\Cinf(L^{-n}\tens\Lambda^\ell TM)
\end{equation*}
associated to the pairing $\W\tens\R\to\W$; this has a rich physical and
geometric interpretation. We have already seen in the previous section that
for $\ell=0$, this is the natural duality between $H_k$ and $H^k$, while the
Leibniz rule for $\ell=1$ shows that $\OpTw_k$ and $\OpTw^k$ are adjoint. For
$\ell=2$, we have a bivector density, integrable over (cooriented) codimension
two submanifolds and the Leibniz rule therefore shows that in the flat case,
$\OpTw_k$ and $\OpTw^{k+1}$ are ``adjoint in codimension one''. In particular,
if $k=0$, $\ell=2$, this is an example of dynamic helicity lowering, using a
twistor to construct a bivector density from a dynamic field. In the flat
case, the Leibniz rule shows that this will be divergence-free wherever the
dynamic field is source-free.  Integrating over a cooriented compact
codimension two submanifold then gives a conserved quantity. Hence on a
simply connected region, the dynamic field and codimension two submanifold
define a real valued linear map on $\W$, i.e., an element of $\W\dual$. This
is one motivation for viewing twistors as ``charges''.

\subsection*{Curvature corrections} If the geometry is not flat, then
some of the above statements only hold up to curvature terms. However, the
curvature corrections appearing in~\ref{compcurv},~\ref{leib} are sometimes
trivial even for non-flat geometries. For instance, when the bilinear pairings
and operators are zero or first order, then there may not be enough
derivatives for curvature to contribute. Also, some parts of the curvature
might not act on certain modules, so that partial flatness assumptions suffice
to kill the curvature terms. Finally, there is the simple observation that if
$\OpTw_k\alpha=0$ then $[K_M]\cptw\alpha=0$, which can be quite a strong
condition if the cup product has low order.

\section{Explicit examples in conformal geometry}

The BGG calculus permits one to carry out many calculations without worrying
too much what the linear and bilinear differential operators are. Nevertheless,
in applications, it is sometimes desirable to determine the operators
explicitly, in terms, for instance, of a chosen Weyl connection, although this
will only be feasible if the operators and pairings have low order. The
Neumann series definition for $Q$, together with the results of
Kostant~\cite{Kostant}, provides one method to carry out these calculations.
However, for many examples, particularly in conformal geometry, one can
proceed more directly, by guessing what the operators and pairings are, using
the considerable collective experience in the literature~\cite{Baston2,
Baston3,Branson,CSS1,CSS3,Thesis,FG,Fegan,Slovak}. The work of~\cite{CSS4} and
this paper contributes to this process by asserting the existence of the BGG
operators, pairings, and Leibniz rules, so that one knows what to look for.

\subsection*{Twistor invariants} Let $\phi$ be a twistor spinor, i.e.,
a solution of the twistor equation for the spin representation of
$\g=\so(V)$. Then $\phi$ is a spinor field with conformal weight $1/2$
satisfying the equation $D_X\phi=\frac1n c(X)\Dirac\phi$ for vector fields
$X$, where $D$ is an arbitrary Weyl connection (the equation is independent of
this choice), $n$ is the dimension of the conformal manifold, $\Dirac=c\circ
D$ is the Dirac operator and $c(X)$ is Clifford multiplication (with
$c(X)^2=\ip{X,X}\iden$). The twistor operator $\OpTw_0$ in this case is given
by the difference $D\phi-\frac1n c(.)\Dirac\phi$; the Lie algebra homology
bundle here is the Cartan product of the cotangent bundle and the spinor
bundle, which is the kernel of Clifford multiplication by $1$-forms.

In four or more dimensions $[K_M]$ is given by the Weyl curvature
$W^\conf\in\Cinf(\Lambda^2T\dual M\cartan\so(TM))$. If $\phi$ is a twistor
spinor then $0=\OpTw_1\OpTw_0\phi=W^\conf\cptw\phi$. This pairing is zero
order: one can check directly that twistor spinors satisfy
$W^\conf_{X,Y}\act\phi=0$ for any vector fields $X,Y$.

The zero order cup products of two solutions $\phi,\psi$ take values in
bundles of conformal weight one. These include
$\omega(\phi,\psi)=\phi\cartan\psi$, which is a (weight $1$)
middle-dimensional form, $X(\phi,\psi)=\ip{c(.)\phi,\psi}$,
which is a vector field, and $\mu(\phi,\psi)=\ip{\phi,\psi}$, which is a
weight $1$ scalar. With our convention for Clifford multiplication these are
all symmetric in $\phi$ and $\psi$. One readily verifies that $\omega$
satisfies a first order twistor equation, that $X$ is a conformal vector
field, and that $\mu$ has vanishing conformal trace-free Hessian (i.e.,
$\sym_0(D^2\mu+r^D\mu)=0$, where $r^D$ is the normalized Ricci tensor of
$D$). Only the last of these has enough derivatives for Weyl curvature to
enter, but it does not act on twistor spinors. One can check
explicitly:
\begin{align*}
D^2_{X,Y}\ip{\phi,\psi}&=
\tfrac1n\bigl(\ip{c(Y)D_X\Dirac\phi,\psi}+\ip{\phi,c(Y)D_X\Dirac\psi}\bigr)\\
&\quad+\tfrac1{n^2}\bigl(\ip{c(Y)\Dirac\phi,c(X)\Dirac\psi}
+\ip{c(X)\Dirac\phi,c(Y)\Dirac\psi}\bigr)\\
&=-\tfrac12\bigl(\ip{c(Y)c(r^D(X))\phi,\psi}+\ip{\phi,c(Y)c(r^D(X))\psi}\bigr)
+\tfrac2{n^2}\ip{X,Y}\ip{\Dirac\phi,\Dirac\psi}\\
&=-r^D(X,Y)\ip{\phi,\psi}+\tfrac2{n^2}\ip{X,Y}\ip{\Dirac\phi,\Dirac\psi}.
\end{align*}
As an application, wherever $\ip{\phi,\psi}$ is nonzero, it defines a
compatible Einstein metric~\cite{BFGK}.

There is also a first order cup product taking values in the trivial
representation, given by
$C(\phi,\psi)=\ip{\Dirac\phi,\psi}-\ip{\phi,\Dirac\psi}$.  One can again check
directly that $dC=0$, as a consequence of the twistor equation. Note that,
with our convention for Clifford multiplication, $\divg^D
X(\phi,\psi)=\ip{\Dirac\phi,\psi}+\ip{\phi,\Dirac\psi}$, so the Dirac operator
is skew adjoint. This convention also makes $C(\phi,\psi)$ skew in $\phi$ and
$\psi$. However, complexifying if necessary, we may assume the spinor bundle
has an orthogonal complex structure $J$, and define
$C(\phi)=C(\phi,J\phi)$. This is the quadratic invariant of Friedrich and
Lichnerowicz~\cite{BFGK,Lich}.

Friedrich also found a quartic invariant (see~\cite{BFGK}). Many similar
invariants can be obtained by iterating the cup product, i.e., for suitable
pairings, $\phi\cptw(\phi\cptw(\phi\cptw\phi))$ will be a nontrivial
scalar. The details are quite complicated, but one obtains a hierarchy of
quartic invariants by considering the cup products factoring through
$k$-forms, for $k=0,1,2,...$.  The first of these is Friedrich's quartic
invariant.

\subsection*{Helicity} Twistor spinors have been systematically exploited to
study massless field equations (see~\cite{PR}). The focus has been mainly on
first order field equations, where zero and first order pairings with twistor
spinors are used to raise or lower the helicity of massless fields. Most of
these pairings are cup products: apart from the helicity $0$ and helicity
$\pm1/2$ equations (given by the conformal Laplacian and Dirac operator
respectively), the massless field equation is the dynamic equation in a BGG
sequence, and the helicity is $\pm(w+1)$ where $w$ is the conformal weight of
the twistor bundle $H_0(W)$.  Helicity $\pm1$ corresponds to the deRham
complex and electromagnetic fields.

\subsection*{Geometrical field theories}
We focus on helicity $\pm2$: this is the expected helicity of linear theories
of gravity, and because of the close links between gravity and geometry, these
differential equations are of particular interest in conformal differential
geometry. There are (at least) three such theories in four (or more)
dimensions: the massless field equations studied by Penrose, the linearization
of Bach's theory of gravity, and a conformal version of a field theory due to
Fierz.  The corresponding representations of $\g=\so(V)$ are $\Lambda^3V$,
$\g=\so(V)=\Lambda^2V$ and $V$. We have given the Lie algebra homologies
in~\ref{confhom} and the BGG sequences begin:
\begin{align}
\label{pen}&\Cinf(L^1\so(TM))&&\to&&\Cinf(L^1T\dual M\cartan\so(TM))&&\to&&
\Cinf(L^1\Lambda^2T\dual M\cartan\so(TM))\\
\label{bach}&\Cinf(TM)&&\to&&\Cinf(\Sym_0
TM)&&\to&&\Cinf(\Lambda^2T\dual M\cartan\so(TM))\\
\label{fierz}&\Cinf(L^1)&&\to&&\Cinf(L^{-1}\Sym_0TM)&&\to&&\Cinf(L^{-1}T\dual M
\cartan\so(TM))
\end{align}
Here $\Sym_0TM=T\dual M\cartan TM$ denotes the symmetric traceless
endomorphisms---note also that $\so(TM)\cong L^2\Lambda^2T\dual M\cong
L^{-2}\Lambda^2TM$.

In Penrose gravity~\eqref{pen}, the operators are both first order, and the
cup product of twistors with Weyl curvature is zero order, giving an
integrability condition for solving the twistor equation. In four dimensions,
the sequence decomposes into selfdual and antiselfdual parts, each part being
a complex iff the Weyl curvature is antiselfdual or selfdual respectively.
This sequence has been used to give a simple proof of the classification of
compact selfdual Einstein metrics of positive scalar curvature~\cite{Besse}.
It also yields a selfduality result for Einstein-Weyl structures on selfdual
$4$-manifolds~\cite{DMJC2}. In Minkowski space, the Weyl tensor of the
Schwarzschild metric defines a natural dynamic field, which may be viewed as a
linearization of the Schwarzschild solution~\cite{Thesis}.

Since Bach gravity~\eqref{bach} is associated to the adjoint representation,
this sequence is the one arising in the study of deformations and moduli
spaces. The twistor operator here is the first order conformal Killing
operator, whose kernel consists of conformal vector fields. It takes values in
$\Sym_0TM\cong L^2S^2_0T\dual M$, the bundle of linearized conformal metrics.
The second operator is the linearized Weyl curvature, taking values in the
bundle of Weyl tensors.  In four dimensions, the adjoint of the linearized Weyl
operator is sometimes called the Bach operator: it can be applied to the Weyl
curvature itself to give the Bach tensor of the conformal structure.

The composite of the conformal Killing operator and the linear Weyl operator
is a first order cup product with the Weyl curvature, which one readily
computes to be a multiple of $\cL_XW^\conf$: obviously a conformal vector
field has to preserve the Weyl curvature. This cup product is associated to
the Lie bracket pairing $\so(V)\tens\so(V)\to\so(V)$. There is also an inner
product pairing, which gives, for example, a $2$-form-valued first order cup
product between vector fields $K$ and Weyl tensors $W$:
\begin{equation*}
(K\cptw W)(X,Y)=(n-2)\ip{W_{X,Y},DK}-\ip{(\delta^DW)_{X,Y},K}.
\end{equation*}

The twistor equation in Fierz gravity~\eqref{fierz} is second order, and is
the conformal tracefree Hessian mentioned above. Its kernel consists of
Einstein (pseudo)gauges, i.e., a nonvanishing solution gives a length scale
for a compatible Einstein metric. The cup product with curvature is a first
order pairing, sometimes called the Cotton-York operator, since it assigns a
Cotton-York tensor to a (pseudo)gauge. This corresponds to the fact that
Einstein metrics have vanishing Cotton-York tensor.

\subsection*{Helicity lowering} A typical example of helicity lowering
occurs in Penrose gravity. Here the natural zero order contraction of a
dynamic Weyl field $W$ and a Penrose twistor 2-form $\omega$ gives a bivector
density and there is a simple Leibniz rule:
\begin{equation*}
\delta\ip{W,\omega}=\ip{\delta W,\omega}+\ip{W,\Twist\omega}.
\end{equation*}
Even within the framework of conformal geometry, the cup and cap products
considerably generalize these ideas. For example, the analogous process in
Fierz gravity requires first order pairings:
\begin{align*}
\delta\bigl(F(D\mu,.,.)-\tfrac12\mu(\delta^D F)(.,.)\bigr)&=
\mu\divg^D(\Sdivg F)-(\Sdivg F)(D\mu,.)\\
&\quad+\ip{\sym_0(D^2+r^D)\mu,F}+\tfrac12\ip{W^\conf,F}\mu
\end{align*}
where $F\in\Cinf((L^{1-n}TM\cartan\so(TM))$, $\mu\in\Cinf(L^1)$, $\Sdivg F$ is
the symmetric divergence in $L^{1-n}\Sym_0TM$, $\delta^DF$ is the skew
divergence in $L^{-1-n}\Lambda^2TM$, and $W^\conf$ is the Weyl curvature.
Higher order pairings rapidly become very complicated, although a few examples
involving second order pairings can be computed explicitly.

\section*{Appendix: semiholonomic jets and Verma modules}

In this appendix we recall the link with semiholonomic Verma modules. On any
Cartan geometry of type $(\g,P)$, iterating the invariant derivative defines
the following.

\begin{prop}\label{jetop} The map sending a section $s$ to $\bigl(s,\InvDer
s,(\InvDer)^2s,\ldots(\InvDer)^ks\bigr)$ takes values in the
subbundle $\cG\cross_P \hat J^k_0\E$ of
$\Dsum_{j=0}^k\bigl((\tens^j\g_M\dual)\tens E\bigr)$, where $\hat J^k_0\E$ is
the set of all $(\phi_0,\phi_1,\ldots \phi_k)$ in
$\Dsum_{j=0}^k\bigl((\tens^j\g\dual)\tens\E\bigr)$ satisfying \textup(for
$1\leq i<j\leq k$\textup) the equations
\begin{align*}
\quad\phi_{j}(\xi_1,\ldots\xi_i,\xi_{i+1},\ldots\xi_j)
-\phi_{j}(\xi_1,\ldots\xi_{i+1},\xi_i,\ldots\xi_j)
&=\phi_{j-1}(\xi_1,\ldots[\xi_i,\xi_{i+1}],\ldots\xi_j)\\
\phi_{i}(\xi_1,\ldots\xi_i)
+\xi_i.\bigl(\phi_{i-1}(\xi_1,\ldots\xi_{i-1})\bigr)&=0
\end{align*}
for all $\xi_1,\ldots\xi_j\in\g$ with $\xi_i\in\p$.
\proofof{prop}
The equations are those given by the vertical triviality and the Ricci
identity, bearing in mind the horizontality of the curvature $\kappa$
(see Proposition~\ref{basic}).
\end{proof}
This map is sometimes called a ``semiholonomic jet operator'', since it
identifies the semiholonomic jet bundle $\hat J^k E$ with the associated
bundle $\cG\cross_P \hat J^k_0\E$. In particular, $\hat J^k_0\E$ is itself the
fibre of the semiholonomic jet bundle at $0=[P]\in G/P$.  This is a minor
variation of the construction of a semiholonomic jet operator given
in~\cite{CSS1,CSS4,ES}, except that we have presented $\hat J^k_0\E$ as a
complicated subspace of an easily defined $\p$-module, whereas in these
references, $\hat J^k_0\E$ is given as a complicated $\p$-module structure on
an easy vector space, namely
$\Dsum_{j=0}^k\bigl(\tens^j(\g/\p)\dual\bigr)\tens\E$.
The equations defining $\hat J^k_0\E$ have a natural algebraic interpretation
in the dual language of semiholonomic Verma modules.

\begin{defn}\tcite{Baston3,ES} Let $\g$ be a Lie algebra with subalgebra $\p$.
\begin{enumerate}
\item The \emphdef{semiholonomic universal enveloping algebra} $U(\g,\p)$ of
$\g$ with respect to $\p$ is defined to be the quotient of the tensor algebra
$\Tens(\g)$ by the ideal generated by
\begin{equation*}
\{\xi\tens\chi-\chi\tens\xi-[\xi,\chi]:\xi\in\p,\chi\in\g\}.
\end{equation*}
We denote by $U^k(\g,\p)$ the filtration given by the image of
$\Dsum_{j=0}^k(\tens^j\g)$, which is compatible with the algebra structure.

Note that $U(\g,\g)$ is the usual universal enveloping algebra $U(\g)$, that
$U(\p)$ is a subalgebra of $U(\g,\p)$, and that $U(\g,\p)$ is a
$U(\g)$-bimodule.
\item Let $\E\dual$ be a $\p$-module.  Then the \emphdef{semiholonomic Verma
module} associated to $\E\dual$ is the $U(\g)$-module given by $\hat
V(\E\dual)= U(\g,\p)\tens_{U(\p)}\E\dual$.  We denote by $\hat V^k(\E\dual)$
the filtration given by the image of $U^k(\p)\tens\E\dual$.
\end{enumerate}
\end{defn}
As a filtered vector space, $\hat V(\E\dual)$ is naturally isomorphic to
$\bigl(\Tens(\g/\p)\bigr)\tens\E\dual$. The induced action of $\g$ on
$\bigl(\Tens(\g/\p)\bigr)\tens\E\dual$ is most easily described by choosing a
complement $\T$ to $\p$ in $\g$ so that $\g/\p$ is isomorphic to $\T$. Then
the action of $\xi\in\g$ on $v_1\tens\cdots\tens v_k\tens
z\in\tens^k\T\tens\E\dual$ is obtained by tensoring on the left with $\xi$,
then commuting the $\p$ component $\xi_{\p}$ past all the $v_j$'s so that it
can act on $z$. This introduces Lie bracket terms $[\xi_{\p},v_j]$, whose
$\p$-components must in turn be commuted to the right. This process is then
repeated until no elements of $\p$ remain in the tensor product.

If we define $\hat J^\infty_0\E$ be the inverse limit with respect to the
natural maps $\hat J^k_0\E\to\hat J^{k-1}_0\E$ then the equations defining
$\hat J^k_0\E$ imply that $\hat J^\infty_0\E$ is the subspace of
$\bigl(\Tens\g\dual\bigr)\tens\E$ such that the pairing with
$\bigl(\Tens\g\bigr)\tens\E\dual$ descends to $\hat V(\E\dual)$. Comparing
dimensions, we see that $\hat V^k(\E\dual)\isom(\hat J^k_0\E)\dual$. This
is why the dual of $\E$ is used in the definition of the Verma module.

The advantage of semiholonomic jets is that they are defined purely in terms
of the $1$-jet functor $J^1$, the natural transformation $J^1E\to E$, and some
abstract nonsense. The construction of the $\Pi$-operators was entirely first
order, and so can be carried out formally on the infinite semiholonomic jet
bundle, rather than on smooth sections. This is perhaps the easiest way to see
that the symbols of the operators are the same for all geometries of a given
type, since the semiholonomic Verma module homomorphisms are the same. One
also sees that the operators of the BGG sequences are all strongly invariant
in the sense that they are defined by semiholonomic Verma module
homomorphisms, and so can be twisted with an arbitrary
$\g$-module~\cite{CSS4}.

In the flat case, the first equation in Proposition~\ref{jetop} holds for
$\xi_i\in\g$ (not just $\p$) and so one can work with the usual holonomic jets
and Verma modules. Working dually with $V(\E\dual)$, instead of $J^\infty E$
or $\Cinf(E)$ as we have done here, leads to a cup coproduct on BGG
resolutions of parabolic Verma modules, as described in~\cite{Thesis}. The
constructions of section~\ref{Ainf} now equip the family of resolutions with
an $A_\infty$-coalgebra structure.

\end{document}